\documentclass[psamsfonts]{amsart}

\usepackage{graphicx}
\usepackage{amsmath}
\usepackage{amsfonts}
\usepackage{amssymb}
\usepackage{url}
\usepackage{diagrams}

\newtheorem{theorem}{Theorem}
\newtheorem{corollary}[theorem]{Corollary}
\newtheorem{definition}[theorem]{Definition}
\newtheorem{example}[theorem]{Example}
\newtheorem{lemma}[theorem]{Lemma}

\newtheorem{proposition}[theorem]{Proposition}
\newtheorem{remark}[theorem]{Remark}

\begin{document}

\title{Unstable attractors in manifolds}

\author{J. J. S\'anchez-Gabites}
\address{Facultad de Matem\'aticas\\ Universidad Complutense de Madrid\\ 28040 Madrid, Spain.}
\email{\url{jajsanch@mat.ucm.es}}
\thanks{The author is supported by MEC, MTM2006-0825.}
\subjclass[2000]{}
\begin{abstract}
Let $M$ be a locally compact metric space endowed with a continuous flow $\varphi : M \times \mathbb{R} \longrightarrow M$. Frequently an attractor $K$ for $\varphi$ exists which is of interest, not only in itself but also the dynamics in its basin of attraction $\mathcal{A}(K)$. In this paper the class of {\sl attractors with no external explosions}, which is intermediate between the well known {\sl stable attractors} and the extremely wild {\sl unstable attractors}, is studied. We are mainly interested in their cohomological properties, as well as in the strong relations which exist between their shape (in the sense of Borsuk) and the topology of the phase space.
\end{abstract}
\maketitle

\section{Introduction}

Let $\varphi : M \times \mathbb{R} \longrightarrow M$ be a continous flow in a locally compact metric phase space $M$ (in our case $M$ will almost always be a topological manifold). Frequently an {\sl attractor} $K$ for $\varphi$ exists which is of interest, and not only in itself but also the dynamics in its basin of attraction $\mathcal{A}(K)$. Such a general situation can be very complicated, because sophisticated dynamics can occur in $\mathcal{A}(K) - K$.

{\sl Stable attractors} have been thoroughly studied; their properties are well known and a number of papers deal with them, for example \cite{bogatyi}, \cite{sanjurjo3}, \cite{segal}, \cite{hastings1}, \cite{kapitanski1}, \cite{sanjurjo1} and \cite{sanjurjo} among others. The present article explores a somewhat bigger class of attractors, that of {\sl attractors with no external explosions}, which can show some mild form of instability. They have already been treated in \cite{athan}, \cite{athan2} (albeit from a different point of view) and \cite{moronsanjurjoyo1}. Thus these papers can be thought of as first steps towards an understanding of more general attractors.

\subsection{A brief outline of the paper} After giving the relevant definitions (essentially due to Athanassopoulos \cite{athan}), we begin Section 2 by establishing some preparatory topological properties of attractors with only internal explosions. This leads to a simple characterization of them in Proposition \ref{caracteriza_inestabilidad1} and a rough understanding of the dynamics in their basin of attraction in Lemma \ref{describe_a_de_k}.

Next (Section 3) we confine ourselves to the case where the phase space $M$ is a topological manifold without boundary and start a local study of attractors with no external explosions. By this we mean an analysis of some cohomological properties of them and, most notably, an appropriate modification of the well-known result that the inclusion of a stable attractor $K$ in its basin of attraction $\mathcal{A}(K)$ is a shape equivalence. It is stated (Theorem \ref{proposicion_simetrico}) in terms of the cohomology of the pair $(\mathcal{A}(K),K)$ and can be used to detect external explosions.

Section 4 concentrates on the case where the phase space is a surface and is included as a manageable example (due to the low dimension of the phase space) which is also easy to visualize. It is shown that the equality $\chi(K) = \chi( \mathcal{A}(K) )$ characterizes, among connected attractors, those with no external explosions (thus a geometrical condition yields a strong dynamical property). A complete description of their shape (in the sense of Borsuk) is given.

Afterwards we change gears and perform a global study of attractors with no external explosions. That is, we explore the interaction between the shape theoretical properties of an attractor and those of the phase space it lives in. It turns out that there exist very strong connections here because of the presence of homoclinic orbits, much in the spirit of Morse--Conley theory (a local dynamical condition, that of being an attractor with no external explosions, is related to the topology of the phase space) but of a somewhat innovative fashion since they cannot be detected by classical tools such as Morse decompositions and Morse equations (we will justify this later on).

In Section 5 we prove that not every manifold can contain an attractor with only internal explosions. A precise statement of this result is Theorem \ref{teorema_cohomologia}, which we follow by a number of corollaries and examples. A partial converse is given in Theorem \ref{proposicion_construye}. Corollary \ref{corolario_cuandoestable} powerfully characterizes stable attractors $K \subseteq \mathbb{R}^2$ as those for which $\chi(K) = \chi(\mathcal{A}(K))$, thus showing that attractors with no external explosions can also be used as tools to study stable ones.

Section 6, which closes this paper, contains a result (Theorem \ref{teorema_esfera}) which further illustrates the strength of the relation between the shape of an attractor $K$ with no external explosions and the phase space. Namely, we prove that if $K$ has the shape of the $n$--sphere then the phase space is homeomorphic either to $\mathbb{S}^n \times \mathbb{S}^1$ or to $\mathbb{S}^n \times_t \mathbb{S}^1$ (the latter being a ``twisted'' product).

\subsection{Requirements} 

\subsubsection{Basic notions about dynamical systems.} Our general reference for the elementary theory of dynamical systems is the book by Bhatia and Szeg\"o \cite{bhatiaszego}, and in particular we shall follow their notation (thus if $\varphi : M \times \mathbb{R} \longrightarrow M$ is a continuous flow on $M$ we abbreviate $\varphi(x,t)$ by $x \cdot t$) save for the limit sets of a point $x \in M$, which we denote $\omega(x)$ and $\alpha(x)$ instead of $\Lambda^{+}(x)$ and $\Lambda^{-}(x)$ (see \cite[3.1 Definition., p. 19]{bhatiaszego}).

The phase space $M$ will be locally compact and metrizable, although most of the paper is set up when $M$ is a manifold.

If $A \subseteq M$ and $x \in A$, we define the {\sl positive prolongational limit set of $x$ relative to $A$} as \[J^{+}(x,A) = \bigcap_{U \in \mathcal{E}_A(x),t \geq 0} \overline{U \cdot [t,+\infty)},\] where $\mathcal{E}_A(x)$ denotes the set of all open neighbourhoods of $x$ in $A$. Clearly $x \in A \subseteq B$ implies $J^{+}(x,A) \subseteq J^{+}(x,B)$ and if $A$ is open in $B$, or more generally a neighbourhood of $x$ in $B$, the equality $J^{+}(x,A) = J^{+}(x,B)$ holds. It is not difficult to check that $\omega(x) = J^{+}(x,\{x\})$ and $J^{+}(x) = J^{+}(x,M)$ (see \cite[Exercises 3.5.1, p. 23 and 4.8.1, p. 28]{bhatiaszego}).

The sets $J^{+}(x,A)$ are closed, invariant, and have the following important property: if $J^{+}(x,A)$ possesses a compact neighbourhood $Q$, then there exist $U \in \mathcal{E}_A(x)$ and $t \geq 0$ such that $U \cdot [t,+\infty) \subseteq Q$. From this it follows without difficulty that $J^{+}(x,A)$ is connected whenever it possesses a compact neighbourhood (in particular, if the phase space is compact). These properties are true in every Hausdorff phase space (not necessarily locally compact or metrizable).

All the definitions above can be dualized and give rise to the {\sl negative prolongational limit set of $x$ relative to $A$}, which we denote $J^{-}(x,A)$. It has the same properties as $J^{+}(x,A)$. A simple proof shows that $y \in J^{-}(x) \Leftrightarrow x \in J^{+}(y)$.

\subsubsection{Isolated invariant sets.} $K$ is an {\sl isolated compact invariant set} (following Conley \cite{conley}) if it possesses a compact neighbourhood $N$ (called an {\sl isolating neighbourhood} for $K$) such that $K$ is the maximal invariant set in $N$, that is $K = \{ x \in N : x \cdot \mathbb{R} \subseteq N \}$. Any isolated set has a neighbourhood basis comprised of isolating neighbourhoods $N$. We shall need to refer to the {\sl asymptotic sets} of $N$, namely $N^{+} = \{ x \in N : x \cdot [0,+\infty) \subseteq N \}$ and $N^{-} = \{ x \in N : x \cdot (-\infty,0] \subseteq N \}$. These are always compact, positively and negatively invariant respectively, and clearly $K = N^{+} \cap N^{-}$.

If $K$ is an isolated compact invariant set, an {\sl isolating block} $N$ for $K$ is an isolating neighbourhood such that $\partial N$ is the union of two compact sets $N^i$ and $N^o$ (called the {\sl entrance} and {\sl exit} sets, respectively) satisfying (1) for every $x \in N^i$ there exists $\varepsilon > 0$ such that $x \cdot [-\varepsilon,0) \subseteq M - N$ and for every $x \in N^o$ there exists $\delta > 0$ such that $x \cdot (0,\delta] \subseteq M - N$, (2) for every $x \in \partial N - N^i$ there exists $\varepsilon > 0$ such that $x \cdot [-\varepsilon,0) \subseteq {\rm int}\ N$ and for every $x \in \partial N - N^o$ there exists $\delta > 0$ such that $x \cdot (0,\delta] \subseteq {\rm int}\ N$. These blocks form a neighbourhood basis of $K$ in $M$ (see \cite{churchill} and \cite{easton}). Moreover, when $M$ is a differentiable $n$--manifold and the flow is also differentiable, there exist isolating blocks $N$ which are $n$--manifolds (with boundary $\partial N$) and such that $N^i, N^o \subseteq \partial N$ are also $(n-1)$--manifolds with common boundary $N^i \cap N^o$ (in $\partial N$).

The {\sl Conley index} of $K$ is defined as the homotopy type $h(K)$ of the pointed quotient $( N/N^o,[N^o] )$, where $N$ is any isolating block for $K$. It can be shown to be independent of the particular $N$ chosen (see, for example, \cite{conley} or \cite{salamon1}). The Conley--Euler characteristic of $K$, to be denoted $\chi h(K)$, is the Euler characteristic of the pair $( N/N^o,[N^o] )$, which is readily seen to agree with $\chi(N,N^o) = \chi(N) - \chi(N^o)$.

\subsubsection{Attractors.} If $K$ is an isolated set, we define its {\sl unstable manifold} as the set $W^u(K) = \{ x \in M : \emptyset \neq \alpha(x) \subseteq K \}$ (this coincides with the definition given in \cite[p. 389]{robinsalamon}). Dually, its {\sl stable manifold} is the set $W^s(K) = \{ x \in M : \emptyset \neq \omega(x) \subseteq K \}$. Clearly both are invariant sets which contain $K$. In case $W^s(K)$ is a neighbourhood of $K$, then $K$ is said to be an {\sl isolated attractor} and its {\sl basin of attraction} is its stable manifold, which we shall usually denote $\mathcal{A}(K)$ (observe that the stable manifold is always defined, but we speak of the basin of attraction only in case $K$ is an attractor). $\mathcal{A}(K)$ is an open neighbourhood of $K$ in $M$ (see \cite[1.8 Theorem., p. 60]{bhatiaszego}).

An isolated attractor $K$ is {\sl stable} if it possesses a neighbourhood basis comprised of positively invariant sets or, equivalently, $J^{+}(x) \subseteq K$ for all $x \in \mathcal{A}(K)$. It is not difficult to show that $K$ is stable if and only if $W^u(K) = K$.

A useful tool is provided by \cite[Theorem 1.25., p. 64]{bhatiaszego}, where it is shown that every attractor $K$ (regardless of its stability) uniquely determines a smallest stable attractor $\widehat{K}$ which contains it and such that $\mathcal{A}(K) = \mathcal{A}(\widehat{K})$. Specifically, $\widehat{K} = \bigcup_{x \in K} J^{+}(x)$, and with this description it can be proved that $\widehat{K}$ is connected if $K$ is. Of course if $K$ is stable, $\widehat{K} = K$.

\subsubsection{Algebraic topology and shape theory.} Some acquaintance with shape theory is convenient, and the reader can find abundant information in \cite{borsukshape2}, \cite{borsukshape1}, \cite{dydak} or \cite{mardesic}. A general knowledge about compact polyhedra (or compact ANR's) is also in order, the books \cite{borsukretractos} and \cite{hu} being a good reference for this.

Our main reference for algebraic topology is the book by Spanier \cite{spanier}. Regarding notation, $\check{H}^{*}$ denotes \v{C}ech cohomology, whereas $H_{*}$ and $H^{*}$ denote singular homology and cohomology, respectively (all unreduced and with coefficients in some commutative ring $R$). Recall that \v{C}ech and singular cohomology agree on ANR's.

If $(X,A)$ is a compact polyhedral pair or, more generally, has the shape of a compact polyhedral pair, then $\check{H}^k(X,A;R)$ is finitely generated for every $k \in \mathbb{N}$ and is nonzero only for finitely many values of $k$ (regardless the coefficient ring $R$). Therefore we can define its {\sl Poincar\'{e} polynomial} $p(X,A;R) = \sum_{k \geq 0} p^k(X,A;R) t^k$, whose $k^{\rm{th}}$ coefficient is \[p^k(X,A;R) = {\rm rk}\ \check{H}^k(X,A;R),\] the rank of the $R$--module $\check{H}^k(X,A;R)$; and its {\sl Euler characteristic} \[\chi(X,A) = \sum_{k \geq 0} (-1)^k {\rm rk}\ \check{H}^k(X,A;\mathbb{Z})\] (see \cite[pp. 172 and 173]{spanier}). By its very definition, $\chi(X,A)$ is obtained evaluating $p(X,A;\mathbb{Z})$ at $t = -1$, and from the universal coefficients theorem it follows that it may also be obtained by evaluating $p(X,A;\mathbb{Z}_2)$ at $t = -1$. We shall make use of the fact, which can be found in \cite[Exercise B.1., p. 205]{spanier}, that $\chi(X,A) + \chi(A) = \chi(X)$ (for the equality to make sense it is enough for any two of the characteristics to be defined, since then the third one is automatically defined too).
\bigskip

The author is grateful to J. M. R. Sanjurjo and M. Castrill\'on L\'opez for their helpful conversations.

\section{Attractors with no external explosions. Definition and basic properties}

Let $K \subseteq M$ be an isolated attractor. We shall call a point $x \in \mathcal{A}(K)$ an {\sl explosion point} if $J^{+}(x) \not\subseteq K$ (let us warn the reader that in \cite{athan} this concept has a different meaning). Stable attractors are those having no explosion points whatsoever. With these preliminaries, it does not seem unreasonable to propose the following definition.

\begin{definition} \label{definicion_solo_internas} (See \cite{athan} and \cite{moronsanjurjoyo1}) An attractor $K$ will be said to have only internal explosions (or no external explosions) if it is isolated and every explosion point in $\mathcal{A}(K)$ belongs to $K$, that is $J^{+}(x) \subseteq K\ \forall\ x \in \mathcal{A}(K) - K$.
\end{definition}

This class of attractors provides us, therefore, with an intermediate stage between the well known picture of stable attractors and the wild one of unstable attractors, and can be used to study both. Let us remark that there do exist unstable attractors which have no external explosions, as the following Example \ref{ejemplo_general} shows.

\begin{example} \label{ejemplo_general} Let $Z$ be a compact topological space and define a flow in $Z \times [-1,1]$ by letting points in $Z \times \{-1\}$ and $Z \times \{1\}$ be fixed and go ``upwards'' from $(z,-1)$ to $(z,1)$ otherwise. Now identify $Z \times \{-1\}$ with $Z \times \{1\}$ in the natural way, that is $(z,-1) \sim (z,1)$, for every $z \in Z$ and call $M$ the resulting quotient space. Then $Z \cong Z \times \{-1\}$ is an unstable attractor with only internal explosions in $M$.
\end{example}

The goal of this section is Lemma \ref{describe_a_de_k}, which gives a fairly complete description of the dynamics in the basin of attraction of an attractor with no external explosions. We need some preparatory results.

\begin{proposition} \label{caracteriza_inestabilidad1} Let $M$ be a locally compact metric space. Assume $K \subseteq M$ is an isolated attractor in $M$. The following statements are equivalent:
\begin{enumerate}
	\item $K$ has only internal explosions.
	\item \label{caracteriza_2} $J^{-}(x) \subseteq K \ \forall\ x \in \widehat{K} - K$.
\end{enumerate}
In case any of these holds, the stabilization of $K$ is $\widehat{K} = W^u(K)$.
\end{proposition}
\begin{proof} (1. $\Rightarrow$ 2.) Let $x \in \widehat{K} - K$, we prove that $J^{-}(x) \subseteq K$. Since $M$ is a locally compact Hausdorff space, its one--point compactification $M_{\infty}$ obtained by adjoining the ideal point $\infty$ to $M$ is a compact Hausdorff space. Morever, the flow in $M$ can be extended to $M_{\infty}$ just by letting $\infty$ be fixed. We shall denote $J^{-}_{\infty}(x)$ the negative prolongational limit set of $x$ in $M_{\infty}$, that is \[J^{-}_{\infty}(x) = \bigcap_{U \in \mathcal{E}_{M_{\infty}}(x),t \geq 0} {\overline{U \cdot (-\infty,-t]}}^{M_{\infty}}.\] Clearly $J^{-}_{\infty}(x) \cap M = J^{-}(x)$, so either $J^{-}_{\infty}(x) = J^{-}(x)$, if the latter is compact, or $J^{-}_{\infty}(x) = J^{-}(x) \cup \{\infty\}$, if not.

$J^{-}_{\infty}(x)$ is connected (because it has a compact neighbourhood, namely $M_{\infty}$) and meets $K$. To see this it suffices to observe that $\alpha(x)$ is a nonempty compact invariant set contained in $\widehat{K}$ (therefore in $\mathcal{A}(K)$) so choosing any $y \in \alpha(x)$ we get $\emptyset \neq \omega(y) \subseteq K \cap \alpha(x) \subseteq K \cap J^{-}_{\infty}(x)$. However $J^{-}_{\infty}(x)$ has empty intersection with $\mathcal{A}(K) - K$. Certainly, if there existed $y \in (\mathcal{A}(K) - K) \cap J^{-}_{\infty}(x)$ then, since $y \in M$, we would have $y \in J^{-}(x)$ so $x \in J^{+}(y)$. But the last set is contained in $K$ by hypothesis, contradicting our choice of $x \not\in K$. Therefore $J^{-}_{\infty}(x)$ is a connected set contained in $\overline{\mathcal{A}(K)}^{M_{\infty}}$ which meets $K$ but not $\mathcal{A}(K) - K$, so necessarily $J^{-}_{\infty}(x) \subseteq K$. In particular $J^{-}(x) = J^{-}_{\infty}(x) \subseteq K$.

(2. $\Rightarrow$ 1.) Let $x \in \mathcal{A}(K) - K$. Clearly $J^{+}(x) \subseteq \widehat{K}$ because $\widehat{K}$ is a stable attractor. However $J^{+}(x)$ cannot meet $\widehat{K} - K$ since otherwise there would exist $y \in \widehat{K} - K$ such that $y \in J^{+}(x)$, so $x \in J^{-}(y)$ which contradicts the hypothesis. Hence $J^{+}(x) \subseteq K$.

Finally we shall derive the equality $\widehat{K} = W^u(K)$ for isolated attractors with only internal explosions from condition \ref{caracteriza_2}. Since $K \subseteq \widehat{K}$, we have $W^u(K) \subseteq W^u(\widehat{K}) = \widehat{K}$, the latter equality being due to the fact that $\widehat{K}$ is stable. For the reverse inclusion, if $x \in \widehat{K}$ then $\emptyset \neq \alpha(x) \subseteq \widehat{K}$ (because $\widehat{K}$ is a compact invariant set) and, moreover, the hypothesis implies $\alpha(x) \subseteq J^{-}(x) \subseteq K$ whenever $x \not\in K$. Therefore $x \in W^u(K)$ if $x \in \widehat{K} - K$ so (since trivially $K \subseteq W^u(K)$) we conclude $\widehat{K} \subseteq W^u(K)$.
 \end{proof}

If we agree to call a point $x \in \mathcal{A}(K)$ (or its orbit) {\sl homoclinic} if $\emptyset \neq \alpha(x), \omega(x) \subseteq K$ (which seems a convenient generalization of the standard notion of a homoclinic orbit), the unstable manifold of $K$ is simply the set of all homoclinic orbits. We shall prefer the notation $\mathcal{H}(K)$ for this set, because we find it to be more intuitive. Proposition \ref{caracteriza_inestabilidad1} above shows that, when $K$ is an attractor with no external explosions, $\widehat{K} = \mathcal{H}(K)$.

Athanassopoulos introduced in \cite{athan2} a classification of isolated attractors $K$ in terms of an ordinal number, called the {\sl instability depth} of the attractor. It measures the complexity of the flow in $\mathcal{A}(K) - K$, and in particular attractors with only internal explosions are precisely the ones with the lowest instability depth (save for stable attractors), though we shall not prove this (it is not difficult). Thus they exhibit the mildest possible instability, as we said above.

\begin{corollary} \label{observacion_esparalelizable} Let $M$ be a locally compact metric space and let $K \subseteq M$ be an isolated attractor with only internal explosions (whether stable or not). Then $\mathcal{A}(K) - K$ and $\widehat{K} - K$ are parallelizable (see \cite[2.1 Definition., p. 48]{bhatiaszego}\/).
\end{corollary}
\begin{proof} Consider the stabilization $\widehat{K}$ of $K$. By \cite[Theorem 2.7., p. 70]{bhatiaszego} its region of attraction $\mathcal{A}(\widehat{K}) = \mathcal{A}(K)$ is a separable set and consequently so is its open subset $\mathcal{A}(K) - K$. Moreover, $\mathcal{A}(K) - K$ is locally compact because it is an open subset of $M$, so by \cite[Theorem 2.6., p. 49 and Theorem 1.8., p. 47]{bhatiaszego} together with the first assertion in Proposition \ref{caracteriza_inestabilidad1} it follows that $\mathcal{A}(K) - K$ is parallelizable.

The proof for $\widehat{K} - K$ is similar. \end{proof}

\begin{remark} \label{obs_yanose} An easy consequence of the fact just proved that $\mathcal{A}(K)$ is parallelizable is that if $\Sigma \subseteq \mathcal{A}(K) - K$ is a compact section of $\mathcal{A}(K) - K$ (that is, the trajectory of every $x \in \mathcal{A}(K) - K$ meets $\Sigma$ exactly in one point $\Sigma \cdot x$), then $\varphi|_{\Sigma \times \mathbb{R}} : \Sigma \times \mathbb{R} \longrightarrow \mathcal{A}(K) - K$ is a homeomorphism (we say then that $\Sigma$ is a {\sl topological section} of $\mathcal{A}(K) - K$). A similar statement holds for $\widehat{K} - K$.
\end{remark}

\begin{lemma} \label{describe_a_de_k} Let $M$ be a locally compact, locally connected metrizable phase space. Assume $K \subseteq M$ is an isolated attractor with only internal explosions. Then:
\begin{enumerate}
	\item $\mathcal{A}(K) - K$ has finitely many connected components $C_1, \ldots, C_s$. These can be arranged so that for some $r \leq s$
	\begin{enumerate}
		\item $\widehat{K} = \mathcal{H}(K) = K \bigcup (C_1 \cup \ldots \cup C_r)$,
		\item $\mathcal{A}(K) = \widehat{K} \bigcup (C_{r+1} \cup \ldots \cup C_s)$.
	\end{enumerate}
	\item There exists a neighbourhood basis of $K$ comprised of isolating neighbourhoods $N$ such that $\partial N$ has finitely many connected components $n_1^{-}, \ldots, n_r^{-} \subseteq N^{-}$ and $n_1^{+}, \ldots, n_r^{+}, \ldots, n_s^{+} \subseteq N^{+}$ which satisfy:
	\begin{enumerate}
		\item For every component $C_j$ of $\mathcal{A}(K) - K$ with $j \leq r$ both $n_j^{-}$ and $n_j^{+}$ are topological sections of $C_j$.
		\item For every component $C_j$ of $\mathcal{A}(K) - K$ with $r < j \leq s$ the set $n_j^{+}$ is a topological section of $C_j$.
	\end{enumerate}
\end{enumerate}
\end{lemma}
\begin{proof} (1) Let $\mathcal{A}(K) - K = \bigcup_{j \in J} C_j$ be the decomposition of $\mathcal{A}(K) - K$ as the disjoint union of its connected components. Choose any isolating neighbourhood $N$ for $K$ contained in $\mathcal{A}(K)$ and observe that, since the compact set $\partial N$ is contained in $\bigcup_{j \in J} C_j$ and every $C_j$ is open because $M$ is locally connected, there exist a finite number of components $C_1, \ldots, C_s$ such that $\partial N \subseteq \bigcup_{j=1}^s C_j$. Now we shall show that every component of $\mathcal{A}(K) - K$ must meet $\partial N$, which will imply that in fact $\mathcal{A}(K) - K = \bigcup_{j=1}^s C_j$. Let $C$ be a component of $\mathcal{A}(K) - K$ and assume that $C \cap \partial N = \emptyset$, so either $C \subseteq N$ or $C \cap N = \emptyset$. The first case is impossible because $C$ is invariant under the flow (being a component of the invariant set $\mathcal{A}(K) - K$) and, since $N$ isolates $K$, we would have $C \subseteq K$. The second one is also ruled out because for any $x \in C$ the set $\omega(x) \subseteq \overline{C}$ would be disjoint from $K$, contradicting the fact that $C$ is contained in $\mathcal{A}(K)$. Therefore $\mathcal{A}(K) - K = \bigcup_{j=1}^s C_j$.

If $x \in \widehat{K} - K$ then $J^{-}(x) \subseteq K$ by Proposition \ref{caracteriza_inestabilidad1} and there exist $t \in \mathbb{R}$ and an open neighbourhood $U$ of $x$ such that $U \cdot (-\infty,-t] \subseteq N$ (since $N$ is a compact neighbourhood of $K$, thus of $J^{-}(x)$). This implies that $\emptyset \neq \alpha(y) \subseteq N$ for every $y \in U$, hence $\emptyset \neq \alpha(y) \subseteq K$ because $N$ isolates $K$. Consequently $U \subseteq W^u(K) = \widehat{K}$ and $\widehat{K} - K$ is open in $\mathcal{A}(K) - K$. Since it is clearly closed in $\mathcal{A}(K) - K$ too, any component of $\mathcal{A}(K) - K$ must be either contained in $\widehat{K}$ or disjoint with it. An adequate relabeling allows us to assume that $C_1, \ldots, C_r$ are the components of $\mathcal{A}(K) - K$ contained in $\widehat{K} = \mathcal{H}(K)$ and $C_{r+1},\ldots,C_s$ the ones disjoint with $\widehat{K}$. Thus $\mathcal{A}(K) - \widehat{K} = \bigcup_{j=r+1}^s C_j$ and $\widehat{K} - K = \mathcal{H}(K) - K = \bigcup_{j=1}^r C_j$.

(2) By \cite[Lemma 3]{moronsanjurjoyo1} $K$ possesses a neighbourhood basis of isolating blocks $N$ with the property that $N = N^{+} \cup N^{-}$ (we refer the reader to \cite{moronsanjurjoyo1} for a detailed discussion of isolating blocks for attractors with only internal explosions). We shall prove that any of them has the features described in the statement of the Lemma. So let $N \subseteq \mathcal{A}(K)$ be such an isolating block contained in $\mathcal{A}(K)$ and observe that its boundary is the disjoint union of the two closed sets $n^{+} = \partial N \cap N^{+}$ and $n^{-} = \partial N \cap N^{-}$. For every $x \in N^i$ there exists $\varepsilon > 0$ such that $x \cdot (-\varepsilon,0) \subseteq M - N$, so $x \not\in N^{-}$ and necessarily $x \in n^{+}$. Thus $N^i \subseteq n^{+}$ and similarly one establishes that $N^o \subseteq n^{-}$. But since $N^i \cup N^o = \partial N = n^{+} \cup n^{-}$ and $n^{+} \cap n^{-} = \emptyset$ we conclude that $N^i = n^{+}$ and $N^o = n^{-}$. This implies that $n^{+}$ is a topological section of $\mathcal{A}(K) - K$, as we proceed to show now. It will be enough to prove that it is a point--set section (see Remark \ref{obs_yanose}) and moreover the fact that $N$ is an isolating block and $n^{+} = N^i$ guarantees that the trajectory of any $x \in \mathcal{A}(K) - K$ intersects $n^{+}$ in at most one point, so it only remains to prove that it indeed meets $n^{+}$. Clearly given any $x \in \mathcal{A}(K) - K$ there must exist some $t \in \mathbb{R}$ such that $x \cdot t \not\in N$ (otherwise $N$ would not isolate $K$) but, since $x \cdot t$ is attracted by $K$, for some other $t^{+}(x) > t$ the relation $x \cdot t^{+}(x) \in \partial N$ must hold. It follows that $x \cdot t^{+}(x) \not\in N^{-}$, because $x \cdot t \not\in N$, so $x \cdot t^{+}(x) \in n^{+}$ as we wanted to prove. The restriction of the flow $\varphi|_{n^{+} \times \mathbb{R}} : n^{+} \times \mathbb{R} \longrightarrow \mathcal{A}(K) - K$ is thus a homeomorphism and $\mathcal{A}(K) - K$ has $s$ connected components $C_1, \ldots, C_s$, so the set $n^{+}$ must also have precisely $s$ components $n^{+}_1, \ldots, n^{+}_s$ which can be labeled in such a way that each $n_j^{+}$ is a section of $C_j$.

Finally, by means of an argument similar to the one used for $n^{+}$ one can show that $n^{-}$ is a topological section for $\widehat{K} - K$. From this the Lemma follows.
 \end{proof}

\begin{remark} \label{obs_global} We keep the notation and hypotheses of Lemma \ref{describe_a_de_k}.
\begin{enumerate}
	\item $K$ is unstable if, and only if, $r \geq 1$ (since otherwise $\widehat{K} = K$).
	\item The components $C_1, \ldots, C_r$ (when $r \geq 1$) are entirely comprised of homoclinic orbits. Thus we shall call them {\sl homoclinic components}. The remaining components $C_{r+1}, \ldots, C_s$ will be said to be {\sl components of uniform attraction}.
	\item If there are no components of uniform attraction then $K$ will be termed {\sl purely unstable}. This is equivalent to having $\mathcal{A}(K) = \widehat{K}$.
	\item If $M$ is connected and $K$ is purely unstable then $K$ is a global attractor (that is, $\mathcal{A}(K) = M$) and the phase space $M$ is compact.
	\begin{proof} $\mathcal{A}(K) = \widehat{K}$ is an open (being a basin of attraction) and closed (being $\widehat{K}$ compact) subset of the connected set $M$, so $M = \mathcal{A}(K) = \widehat{K}$. In particular $M$ is compact.
	 \end{proof}
	\item Conversely, if $M$ is compact then any global connected isolated attractor with only internal explosions $\widehat{K}$ is purely unstable.
	\begin{proof} Its stabilization $\widehat{K}$ is a stable attractor whose basin of attraction is compact, which requires $\widehat{K} = M$ and implies that $K$ is purely unstable.
	 \end{proof}
\end{enumerate}
\end{remark}

Isolating neighbourhoods satisfying the conditions in part (2) of Lemma \ref{describe_a_de_k} will be called {\sl regular}. Note however that this is not the standard usage of the word.

\section{Attractors with no external explosions: cohomological properties and relations with their basins of attraction}

After having given a rough description of the basin of attraction of an attractor with no external explosions, we concentrate on the attractors themselves. In \cite{moronsanjurjoyo1} it is proved that such an attractor has the shape of a finite polyhedron. Now we want to obtain some properties about their \v{C}ech cohomology and the \v{C}ech cohomology of their basins of attraction.

An {\sl $n$--manifold} will mean a topological manifold of dimension $n$ without boundary. Recall that every such manifold is locally connected, locally compact and, when connected, metrizable (thus we are under the conditions of the previous section).

We begin with a duality result which allows us to treat regular isolating neighbourhoods as if they were manifolds with boundary, at least at homology level. It can be proved either resorting to the theory of generalized manifolds (about which the reader can find information in \cite{bredon1}, \cite{raymond} and \cite{wilder}) or in a more classical fashion, as follows.

\begin{lemma} \label{lema_dualidad} Let $M$ be an $R$--orientable $n$--manifold without boundary. Assume that $K \subseteq M$ is a connected isolated unstable attractor with only internal explosions. If $N$ is a connected regular isolating neighbourhood for $K$ contained in $\mathcal{A}(K)$, then:
\begin{enumerate}
	\item $H^k(N,\partial N;R) = H_{n-k}(N;R)$ (Lefschetz duality for $N$),
	\item $H^k(n^{-};R) = H_{n-1-k}(n^{-};R)$ (Poincar\'e duality for $n^{-}$).
\end{enumerate}
\end{lemma}
\begin{proof} $\mathcal{A}(K) - K$ is a manifold (because it is open in $M$) and therefore an ANR. By Lemma \ref{describe_a_de_k} $n^{+}$ is a topological section (hence a retract) of $\mathcal{A}(K) - K$, and consequently it is an ANR. Similary one shows that $n^{-}$ is an ANR ($n^{-}$ is a retract of $\widehat{K} - K$, which is open in $\mathcal{A}(K) - K$) so $\partial N = n^{+} \cup n^{-}$ is an ANR. In the same fashion it can be proved that $N$ is an ANR, and all this implies that \v{C}ech and singular cohomology agree on $N$, $\partial N$ and the pair $(N,\partial N)$.

Coefficients are taken in $R$.

(1) By Alexander duality on $\mathcal{A}(K)$ we have \[H^k(N,\partial N) = H_{n-k}(\mathcal{A}(K) - \partial N, \mathcal{A}(K) - N).\] Now $\mathcal{A}(K)$ is an open neighbourhood of the compact set $N$, so $\mathcal{A}(K) - \partial N = {\rm int}(N) \bigcup \mathcal{A}(K) - N$, where the union is disjoint and both sets are open. Thus \[H_{n-k}(\mathcal{A}(K) - \partial N, \mathcal{A}(K) - N) = H_{n-k}({\rm int}(N)).\] Now it follows from the proofs in \cite{moronsanjurjoyo1} that the last group is isomorphic to $H_{n-k}(N)$.

(2) Let $C = \widehat{K} - K = n^{-} \cdot \mathbb{R}$, which is an $R$--orientable $n$--manifold because it is open in $\mathcal{A}(K)$ by Lemma \ref{describe_a_de_k}. By Alexander duality \[H_k(C,C - n^{-}) = H^{n-k}(n^{-}).\] The pair $(C,n^{-})$ is homeomorphic to $(n^{-} \times (-1,1),n^{-} \times \{0\})$ since $n^{-}$ is a section of $C$ (Lemma \ref{describe_a_de_k}). Consequently \[H_k(C,C-n^{-}) = H_k(n^{-} \times [-1,1],n^{-} \times \{-1,1\}),\] and the latter group is isomorphic to $H_{k-1}(n^{-})$ by Lemma 5. With the duality relation established above we get $H_k(n^{-}) = H^{(n-1)-k}(n^{-})$. \end{proof}

\begin{proposition} \label{proposicion_hiper} Let $M$ be an $R$--orientable $n$--manifold without boundary. Assume that $K \subseteq M$ is a connected isolated attractor with no external explosions. Then:
\begin{enumerate}
	\item $\check{H}^k(K;R) = 0$ for $k \geq n$,
	\item $r \leq s \leq {\rm rk}\ \check{H}^{n-1}(K;R)$,
\end{enumerate}
where $s$ is the total number of components in $\mathcal{A}(K) - K$ and $r$ is the number of homoclinic components in $\mathcal{A}(K) - K$.
\end{proposition}
\begin{proof} Take coefficients in $R$.

(1) By Alexander's duality $\check{H}^k(K) = H_{n-k}(\mathcal{A}(K),\mathcal{A}(K)-K)$. For $k > n$ we trivially have $H_{n-k}(M,M - K) = 0$. For $k = n$, recalling the fact that $\mathcal{A}(K)$ is connected (because $K$ is) it follows that $H_0(\mathcal{A}(K),\mathcal{A}(K)-K) = 0$.

(2) We prove that $s \leq {\rm rk}\ \check{H}^{n-1}(K)$. Letting $N$ be a connected isolating block for $K$, recall that by Lemma \ref{describe_a_de_k} $r$ is the number of components of $n^{-}$ and $s$ is the number of components of $n^{+}$. Consider the long exact sequence in \v{C}ech cohomology for the triple $(N,\partial N,n^{-})$. Since $N$ is connected, $\check{H}^0(N,n^{-}) = 0$ so \[\ldots \longleftarrow \check{H}^1(N,\partial N) \longleftarrow \check{H}^0(\partial N,n^{-}) \longleftarrow \check{H}^0(N,n^{-}) = 0 \longleftarrow \ldots\] implies that ${\rm rk}\ \check{H}^0(\partial N,n^{-}) \leq {\rm rk}\ \check{H}^1(N,\partial N)$. Now $\partial N = n^{+} \cup n^{-}$ and therefore $\check{H}^0(\partial N,n^{-}) = \check{H}^0(n^{+})$ has rank $s$. By Lefschetz duality on $N$ (Lemma \ref{lema_dualidad}) $H^1(N,\partial N) = H_{n-1}(N)$, so $s \leq {\rm rk}\ H_{n-1}(N)$. Using the universal coefficients theorem ${\rm rk}\ H_{n-1}(N) = {\rm rk}\ H^{n-1}(N)$ and since $N$ is an ANR the latter coincides with $\check{H}^{n-1}(N)$. Finally, the inclusion $K \hookrightarrow N$ is a shape equivalence (\cite[Proposition 5]{moronsanjurjoyo1}), whence $s \leq {\rm rk}\ \check{H}^{n-1}(K)$.
 \end{proof}

Observe in particular that, if $K$ is unstable (thus $r \geq 1$), then ${\rm rk}\ \check{H}^{n-1}(K) \geq 1$ so ${\rm dim}(K) \geq n-1$. This was proved in \cite[Theorem 4.5, p. 166]{athan2} under the assumption that the flow is smooth.

\begin{corollary} \label{corolario_global} Let $M$ be a connected $R$--orientable $n$--manifold without boundary. Assume that $K \subseteq M$ is a connected isolated unstable attractor with no external explosions. If ${\rm rk}\ \check{H}^{n-1}(K;R) = 1$, then $K$ is a global attractor. Further, $M$ is compact and $r = 1$.
\end{corollary}
\begin{proof} By Proposition \ref{proposicion_hiper} we have $1 \leq r \leq s \leq {\rm rk}\ \check{H}^{n-1}(K;R)$ so $r = s = 1$. There exists just one homoclinic component because $r = 1$ but, moreover, $r = s$ implies by Lemma \ref{describe_a_de_k} that $\mathcal{A}(K) = \widehat{K}$. By Remark \ref{obs_global} the attractor $K$ is global and $M$ is compact, because it is connected.
 \end{proof}

It is known that for a stable attractor $K$ the inclusion $K \hookrightarrow \mathcal{A}(K)$ is a shape equivalence (\cite{bogatyi}, \cite{sanjurjo3}, \cite{segal}, \cite{hastings1}, \cite{kapitanski1}, \cite{sanjurjo1}, \cite{sanjurjo}) so the polynomial $p(\mathcal{A}(K),K)$ is trivial. When $K$, more generally, has only internal explosions, the following result holds.
 
\begin{theorem} \label{proposicion_simetrico} Let $M$ be an $R$--orientable $n$--manifold without boundary. Assume that $K \subseteq M$ is a connected isolated attractor with only internal explosions. Then the polynomial $p(\mathcal{A}(K), K ; R )$ has the symmetric structure \[p( \mathcal{A}(K), K ; R ) = a_n t^n + a_{n-1} t^{n-1} + \ldots + a_{n-1} t^2 + a_n t,\] with all the coefficients finite. Moreover $a_n = r$, where $r$ is as in Lemma \ref{describe_a_de_k} the number of homoclinic components of $\mathcal{A}(K) - K$.
\end{theorem}
\begin{proof} Let $N$ be a connected regular isolating block for $K$, so the inclusion $K \hookrightarrow N$ is a shape equivalence. All coefficients are taken in $R$.

Since the inclusion $(\widehat{K},K) \hookrightarrow (\mathcal{A}(K),K)$ is a shape equivalence it follows that $p(\mathcal{A}(K),K) = p(\widehat{K},K)$. By Lemma \ref{describe_a_de_k} we have $\widehat{K} - K = n^{-} \cdot \mathbb{R}$ which is homeomorphic to $n^{-} \cdot (-1,1) = n^{-} \cdot [-1,1] - n^{-} \cdot \{-1,1\}$ and using Lemma \ref{observacion_cociente} we see that \[p(\mathcal{A}(K),K) = p(n^{-} \cdot [-1,1],n^{-} \cdot \{-1,1\}).\] Applying Lemma \ref{lema_suspension} the latter polynomial equals $t p(n^{-})$, which has finite coefficients because $n^{-}$ is an ANR.

By the universal coefficients theorem the free part of $H^k(n^{-})$ is isomorphic to the free part of $H_k(n^{-})$. This implies that both have the same rank $p^k(n^{-})$, but by the duality relation established in Lemma \ref{lema_dualidad}.(2) this rank also coincides with that of $H^{(n-1)-k}(n^{-})$, that is $p^{(n-1)-k}(n^{-})$. Therefore $p(n^{-})$ is symmetric and using the relation $p( \mathcal{A}(K),K ) = t p(n^{-})$ established above it follows that $p( \mathcal{A}(K),K )$ is also symmetric.

Due to this symmetry, the leading coefficient of $p(\mathcal{A}(K),K)$ equals $p^0(n^{-})$ which is the number of connected components of $n^{-}$ or equivalently (by Lemma \ref{describe_a_de_k}) the number of homoclinic components in $\mathcal{A}(K) - K$.
 \end{proof}

Hence we see that the polynomial $p(\mathcal{A}(K),K)$ can be either trivial (if $K$ is stable) or have degree $n = {\rm dim}(M)$ if $K$ is unstable. It provides another method to detect external explosions, as presented in the following Example \ref{ejemplo_eneltoro} (compare also with Theorem \ref{teorema_esfera}):

\begin{example} \label{ejemplo_eneltoro} Let $K \subseteq \mathbb{T}^3$ be an isolated unstable attractor with the shape of $\mathbb{S}^2$. Then $K$ must have external explosions.

We shall calculate $p(\mathcal{A}(K),K)$ and observe that it does not satisfy the symmetry conditions required by Theorem \ref{proposicion_simetrico} (all coefficients are taken in $\mathbb{Z}$). Since $\check{H}^2(K) = \mathbb{Z}$, we have $r = s = 1$ so that $K$ is purely unstable and hence a global attractor (Corollary \ref{corolario_global}). Thus $p(\mathcal{A}(K),K) = p(\mathbb{T}^3,K) = a_3t^3 + a_2t^2 + a_1t$ for some nonnegative integers $a_1,a_2,a_3$. We already know that $a_3 = r = 1$ (Theorem \ref{proposicion_simetrico}). From the exact sequence for the pair $(\mathbb{T}^3,K)$ and using $\check{H}^1(K) = 0$ it follows that $\check{H}^1(\mathbb{T}^3,K) = \check{H}^1(\mathbb{T}^3) = \mathbb{Z}^3$, so $p(\mathcal{A}(K),K) = t^3 + a_2t^2 + 3t$ which contradicts the symmetry required by Theorem \ref{proposicion_simetrico}.
\end{example}

The following Corollary will prove useful later on.

\begin{corollary} \label{corolario_caracteristica} Let $M$ be a manifold without boundary of even dimension $n$. Assume that $K \subseteq M$ is a connected isolated attractor with no external explosions. Then $\chi(\mathcal{A}(K)) = \chi(K)$.
\end{corollary}
\begin{proof} Let $n$ be even and take coefficients $R = \mathbb{Z}_2$, so that $M$ is indeed $R$--orientable. The symmetry of $p(\mathcal{A}(K),K)$ proved in Theorem \ref{proposicion_simetrico} implies that $\chi(\mathcal{A}(K),K) = 0$, so $\chi(\mathcal{A}(K)) = \chi(K) + \chi( \mathcal{A}(K), K ) = \chi(K)$.
 \end{proof}

It is not difficult to see that, if the dimension of the phase space is odd, Corollary \ref{corolario_caracteristica} need not be true any more.
\bigskip

In \cite[Theorem 2., p. 327]{segal} it is proved that a finite dimensional compactum $K$ with the shape of a compact polyhedron can be embedded in euclidean space in such a manner that $K$ is a stable attractor for some suitable flow. It is reasonable to ask whether the same is true if we want $K$ to be unstable and have no external explosions (anticipating things to come, embeddings in manifolds other than $\mathbb{R}^n$ should be allowed, because of Example \ref{ejemplo_esfera}). To finish this section we answer this question in the negative with Example \ref{nopuedeser}. Thus one obtains the somewhat surprising result that there exist compacta which, owing only to their shape, must have external explosions when embedded as unstable attractors in manifolds.

\begin{example} \label{nopuedeser} Consider the wedge sum $\mathbb{S}^2 \vee \mathbb{S}^1 \subseteq \mathbb{R}^3$ and thicken it a little to obtain a compact $3$--manifold $K$ (with boundary) which collapses onto $\mathbb{S}^2 \vee \mathbb{S}^1$ ($K$ is a regular neighbourhood of $\mathbb{S}^2 \vee \mathbb{S}^1$ in $\mathbb{R}^3$). We assert that $K$ cannot be embedded as an unstable attractor with no external explosions in any manifold $M$ without boundary.
\end{example}
\begin{proof} Let us proceed to prove the claim by contradiction. We shall assume that $M$ is connected (this is no loss of generality, because $K$ is connected). Taking coefficients in $\mathbb{Z}_2$, it is straightforward to check that $\check{H}^0(K) = \check{H}^1(K) = \check{H}^2(K) = \mathbb{Z}_2$, the rest of the groups being null. If we let $S, L \subseteq K$ be copies of $\mathbb{S}^2, \mathbb{S}^1 \subseteq \mathbb{S}^2 \vee \mathbb{S}^1$ slightly displaced so as to make them disjoint, is clear that the inclusion induced homomorphisms $H^2(K) \longrightarrow H^2(S)$ and $H^1(K) \longrightarrow H^1(L)$ are isomorphisms, and similarly for homology. Thus $S$ and $L$ are geometric representatives of the (co)homology of $K$.

Denote $m = {\rm dim}(M)$. Then $m \geq 3$ since $\check{H}^k(K)$ must be zero for $k \geq m$ (by Proposition \ref{proposicion_hiper}) and, on the other hand, $m \leq 3$ since $\check{H}^{m-1}(K)$ must be nonzero (by the same token). Thus $m = 3$. Corollary \ref{corolario_global} asserts that $K$ is a global attractor because ${\rm rk}\ \check{H}^2(K) = 1$. Further $M$ is compact and $r = 1$.

Now we compute $p(M,K)$. By Theorem \ref{proposicion_simetrico} we see that $p(M,K) = t^3 + at^2 + t$ for some $a \geq 0$. Evaluating this equality at $t = -1$ yields $\chi(M,K) = -2 + a$. It is a well known consequence of Poincar\'e duality that $\chi(M) = 0$ because $M$ is a compact $3$--manifold without boundary, and direct calculation shows that $\chi(K) = 1$. Hence $-2 + a = \chi(M,K) = \chi(M) - \chi(K) = -1$ and consequently $a = 1$, so $p(M,K) = t^3 + t^2 + t$.

Let $N$ be an isolating block for $K$ contained in $\mathcal{A}(K)$. It can be proved that $N$ is a compact $3$--manifold (the argument is given in Step 1 of Theorem \ref{teorema_esfera}). Its boundary $\partial N$ is the disjoint union of the closed surfaces $n^{-}$ and $n^{+}$, which are homeomorphic because both are sections of $M - K$ (see Lemma \ref{describe_a_de_k}). Further $p(M,K) = tp(n^{-})$ (recall the proof of Theorem \ref{proposicion_simetrico}) so $p(n^{-}) = t^2 + t + 1$ and (since the only surface which has this Poincar\'{e} polynomial is $\mathbb{RP}^2$, the projective plane) we conclude that $\partial N$ is the disjoint union of two projective planes.

$S \cong \mathbb{S}^2$ possesses a small open product neighbourhood $U \cong \mathbb{S}^2 \times (-1,1) \subseteq {\rm int}(K) \subseteq {\rm int}(N)$, whose boundary $\partial U$ is the disjoint union of two copies of $\mathbb{S}^2$. We will show that the compact $3$--manifold $N - U$ (with boundary $\partial N \cup \partial U$) has two connected components, each of which has a boundary consisting of a $2$--sphere (coming from $\partial U$) and a projective plane $\mathbb{RP}^2$ (coming from $\partial N$). However this is not possible, because the boundary of an odd--dimensional manifold (such as one of the components of $N - U$) must have even Euler characteristic (this is easy to prove) but $\chi(\mathbb{S}^2) + \chi(\mathbb{RP}^2) = 3$. Thus a contradiction will be obtained so the hypothesis that $K$ is unstable but has no external explosions is untenable.
\medskip

{\it Claim.} There is an isomorphism $H_1(N,\partial N) \longrightarrow H_1(N,N-S)$ induced by inclusion.

\begin{proof} Let $N'$ be $N$ with an open external collar attached on $\partial N$. Both inclusions $K \hookrightarrow N$ and $N \hookrightarrow N'$ are homotopy equivalences, so the inclusion induced homomorphism $H^2(N') \longrightarrow H^2(S)$ is an isomorphism. By Alexander duality (and its naturality properties) we have a commutative diagram
\begin{diagram}
H^2(N') & \rTo^{\cong} & H^2(S) = \mathbb{Z}_2 \\
\dTo <{\cong} & & \dTo <{\cong} \\
H_1(N',N'-N) & \rTo & H_1(N',N'-S)
\end{diagram}
which implies that the lower arrow is also an isomorphism. Now $(N',N'-N)$ has the same homotopy type as $(N',N'-{\rm int}(N))$, and in turn the latter deformation retracts onto $(N,\partial N)$. Since $(N',N'-S)$ is carried onto $(N,N-S)$ under this deformation, we conclude that the inclusion induced homomorphism $H_1(N,\partial N) \longrightarrow H_1(N,N-S)$ is an isomorphism.
\end{proof}
\medskip

{\it Claim.} Every component of $N - U$ meets $\partial N$.

\begin{proof} From the exact sequence for the triple $(N,N-S,\partial N)$
\begin{multline*} \ldots \longrightarrow H_1(N,\partial N) \stackrel{\cong}{\longrightarrow} H_1(N,N-S) \longrightarrow \\ \longrightarrow H_0(N-S,\partial N) \longrightarrow H_0(N,\partial N) = 0 \end{multline*}
it follows that $H_0(N-S,\partial N) = 0$. Therefore every component of $N - S$ meets $\partial N$, and so does every component of $N - U$.
\end{proof}
\medskip

{\it Claim.} $N - U$ has precisely two connected components.

\begin{proof} Since $N - U$ deformation retracts onto $N - S$, it will be enough to show that the latter has two connected components. The chain of inclusions $L \subseteq N - S \subseteq N$ gives rise, in homology, to the composition $H_1(L) \longrightarrow H_1(N-S) \longrightarrow H_1(N)$ which we know to be an isomorphism by the choice of $L$. Hence $H_1(N - S) \longrightarrow H_1(N)$ is surjective and from the exact sequence
\begin{multline*} \ldots \longrightarrow H_1(N-S) \longrightarrow H_1(N) \longrightarrow H_1(N,N-S) \longrightarrow \\ \longrightarrow \widetilde{H}_0(N-S) \longrightarrow \widetilde{H}_0(N) = 0 \end{multline*}
it follows that $\widetilde{H}_0(N-S)$ is isomorphic to $H_1(N,N-S) = \mathbb{Z}_2$, thus establishing that $N - S$ has two connected components.
\end{proof}

Consequently $N - U$ has two connected components, each of which contains one component of $\partial N$ (one projective plane) and one of $\partial U$ (a $2$--sphere). This yields the contradiction we stated above and finishes the proof.  \end{proof}

\section{Attractors with no external explosions in surfaces}

We devote this section to the case where the phase space $M$ is a surface (not necessarily compact). By Corollary \ref{corolario_caracteristica} we already know that if $K \subseteq M$ is a connected attractor with no external explosions then $\chi(K) = \chi(\mathcal{A}(K))$. Our first result shows that this simple condition provides, in fact, a characterization of these attractors.

\begin{theorem} \label{teorema_caracteriza} Let $M$ be a $2$--manifold without boundary and let $K \subseteq M$ be a connected isolated attractor. Then $K$ has only internal explosions if, and only if, $\chi(K) = \chi( \mathcal{A}(K) )$.
\end{theorem}
\begin{proof} We only need prove that $\chi(K) = \chi(\mathcal{A}(K))$ implies that $K$ has no external explosions. Assume (a justification for this is provided later) that the flow is topologically conjugate to one of class $\mathcal{C}^1$, at least in some open neighbourhood of $K$.
\medskip

{\it Step 1.}  Let $N$ be a connected isolating block for $K$. Denote $n^{+} = N^{+} \cap \partial N$ and $n^{-} = N^{-} \cap \partial N$. Then every component of $\partial N$ meets $n^{-}$ or $n^{+}$.
\begin{proof} Consider the following commutative diagram, which is built up from the exact sequences for the pairs $(N^o,n^{-})$, $(N^{-},n^{-})$ and $(N,N^o)$.
\begin{diagram}
& & 0 & & & & \\
& & \dTo & & & & \\
& & \check{H}^0(N^o,n^{-}) & & & & \\
& & \dTo & & & & \\
\check{H}^1(N,N^o) & \lTo & \check{H}^0(N^o) & \lTo & \check{H}^0(N) & \lTo & \check{H}^0(N,N^o) \\
\dTo <{(1)} & & \dTo <{(2)} & & \dTo <{(3)} & & \dTo <{(4)} \\
\check{H}^1(N^{-},n^{-}) & \lTo & \check{H}^0(n^{-}) & \lTo & \check{H}^0(N^{-}) & \lTo & \check{H}^0(N^{-},n^{-})
\end{diagram}
Arrows (1) and (4) are isomorphisms by \cite[Corollary 2.(a), p. 1437]{sanjurjo0}, and so is (3) because $N^{-}$ is a connected (since $K \hookrightarrow N^{-}$ is a shape equivalence) subset of the connected set $N$. An easy diagram chasing shows then that (2) is a monomorphism and consequently $\check{H}^0(N^o,n^{-}) = 0$. That is, every component of $N^o$ meets $n^{-}$. Similarly one proves that every component of $N^i$ meets $n^{+}$.
\end{proof}
\medskip

{\it Step 2.} $\chi h(K) = \chi(K)$.
\begin{proof} As $\widehat{K}$ (the stabilization of $K$) is a stable attractor with region of attraction $\mathcal{A}(K)$, its Conley--Euler characteristic is $\chi h(\widehat{K}) = \chi(\widehat{K})$ (either take a positively invariant neighbourhood of $\widehat{K}$ as an isolating block to compute its Conley index or see \cite[Corollary 2.(b), p. 1437]{sanjurjo0}). Since the inclusion $\widehat{K} \hookrightarrow \mathcal{A}(\widehat{K}) = \mathcal{A}(K)$ is a shape equivalence, $\chi h(\widehat{K}) = \chi(\widehat{K}) = \chi(\mathcal{A}(K))$.

Every fixed point of $\widehat{K}$ (if any) must lie in $K$, because the latter is an attractor and $\widehat{K}$ is contained in its basin of attraction. Hence by \cite[Corollary.(ii), p. 858]{mccord}, $\chi h(K) = \chi h(\widehat{K}) = \chi(\mathcal{A}(K))$, and using our hypothesis that $\chi(K) = \chi(\mathcal{A}(K))$ we get $\chi h(K) = \chi(K)$.

A warning may be in order here. McCord's paper \cite{mccord} is set up for $\mathcal{C}^1$ flows in compact manifolds. However it is not difficult to see that only differentiability near the invariant set is used. Further, since Conley's index and fixed points are preserved by topological conjugation, the assertions of \cite[Corollary., p. 858]{mccord} are also true for flows which are just topologically conjugate to a class $\mathcal{C}^1$ one near the invariant set. We shall see later on that this is our case.
\end{proof}
\medskip

{\it Step 3.} $K$ has an isolating block of the form $N = N^{+} \bigcup N^{-}$.
\begin{proof} Let $N$ be a connected isolating block for $K$. Then $\chi h(K) = \chi(N,N^o) = \chi(N^{-},n^{-})$, where the second equality follows again from \cite[Corollary 2.(a), p. 1437]{sanjurjo0}, referred to in Step 1. above. Moreover, $\chi(N^{-}) = \chi(K)$ because the inclusion $K \hookrightarrow N^{-}$ is a shape equivalence, so $\chi h(K) = \chi(N^{-},n^{-}) = \chi(N^{-}) - \chi(n^{-}) = \chi(K) - \chi(n^{-})$ and consequently $\chi(n^{-}) = 0$ (recall that $\chi h(K) = \chi(K)$ by Step 2.). Similarly one proves that $\chi(n^{+}) = 0$.

The argument above applies to any connected isolating block for $K$. Since the flow is topologically conjugate to a class $\mathcal{C}^1$ one near $\widehat{K}$ (in particular near $K$) and isolated sets for differentiable flows have isolating blocks which are differentiable manifolds (\cite{easton}), $K$ has isolating blocks which are topological manifolds. Thus we may assume that $N$ is a compact $2$--manifold with boundary $\partial N$, which is a union of circumferences $C_1, \ldots, C_s$. Let $P = n^{+} \bigcup n^{-}$. Then $P$ is a compact subset of $\partial N$ and, by Step 1., every circumference $C_i$ meets $P$. Moreover $\chi(P) = 0$ because it is the disjoint union of two sets with Euler characteristic zero. Now $0 = \chi(P) = \sum_{i=1}^s \chi(P \cap C_i)$ and, since the Euler characterstic of a proper subset of a circumference is positive, it follows that $P \cap C_i = C_i$ for every $1 \leq i \leq s$, so $\partial N = P = n^{-} \bigcup n^{+}$. Thus $N = N^{+} \bigcup N^{-}$.
\end{proof}
\medskip

{\it Step 4.} $K$ has only internal explosions.
\begin{proof} This is easy. Let $N = N^{+} \bigcup N^{-}$ be an isolating block for $K$, as the one obtained above. We have to show that $J^{+}(x) \subseteq K$ for every $x \in \mathcal{A}(K) - K$. Given $x \in \mathcal{A}(K)$, as $\omega(x) \subseteq K$ (because $K$ is an attractor) there exists $t \geq 0$ such that $x \cdot [t,+\infty) \subseteq N$, or $x \cdot t \in N^{+}$. Now $N^{+}$ is a neighbourhood of $x \cdot t$ in $N$ (because $N = N^{+} \bigcup N^{-}$), hence in $M$. Therefore $J^{+}(x) = J^{+}(x,N^{+})$. Finally, since $K$ is a stable attractor for the semiflow obtained by restricting $\varphi$ to $N^{+}$, we have the inclusion $J^{+}(x,N^{+}) \subseteq K$, which finishes the proof.
\end{proof}
\medskip

It only remains to justify the assumption that the flow is topologically conjugate to one which is of class $\mathcal{C}^1$ near $\widehat{K}$. Let $U$ be an open neighbourhood of $\widehat{K}$ in $M$ such that $\overline{U}$ is a compact $2$--manifold with boundary. By a theorem of Beck (\cite{beck}) we can modify the flow in $M$ so that $\partial U$ is comprised of fixed points. Now $U$ is an invariant set, and the trajectory segments contained in it are unaltered (save a reparametrization) by this process. Pasting disks onto $\partial U$ to get rid of the boundary components of $\overline{U}$ (if any) and letting those disks consist of fixed points, a compact $2$--manifold is obtained which carries a flow containing $U$. Using Guti\'{e}rrez's smoothing theorem \cite{gutierrez}, the flow in $U$ is topologically conjugate via a homeomorphism $h$ to a class $\mathcal{C}^1$ flow.
 \end{proof}

Finally, Theorem \ref{teorema_superficies} gives a complete description of global (unstable) attractors, with only internal explosions, in compact surfaces.

\begin{theorem} \label{teorema_superficies} Let $M$ be a compact, connected $2$--manifold without boundary. If $K \subseteq M$ is a connected isolated unstable global attractor with no external explosions, then $\chi(M) \leq 0$ (thus $M$ is neither the sphere nor the real projective plane) and $K$ has the shape of a bouquet of $1 - \chi(M)$ circumferences.
\end{theorem}
\begin{proof} Let us show in the first place that $K$ has the shape of a bouquet of circumferences. Since $K$ is a global attractor and $M$ is compact, by Remark \ref{obs_global} it follows that $K$ is purely unstable so $\mathcal{A}(K) - K$ consists only of homoclinic orbits. Therefore if $N = N^{+} \cup N^{-}$ is an isolating block for $K$ we know from Lemma \ref{describe_a_de_k} that $K$ has the same shape as $M - n^{-}$. Since $n^{-} \times \mathbb{R}$ is homeomorphic to $M - K$, which is a generalized $2$--manifold (because it is a topological $2$--manifold), $n^{-}$ is a generalized $1$--manifold itself (see \cite[Theorem 6., p. 17]{raymond}, which applies only to locally orientable generalized manifolds, and \cite{bredon1}, showing that every generalized manifold is locally orientable) and therefore a topological $1$--manifold by \cite[Theorem 1.2, p. 271]{wilder}. Hence $n^{-}$ is the disjoint union of $r \geq 1$ circumferences (alternatively one can use Gutierrez's smoothing theorem \cite{gutierrez} as in Theorem \ref{teorema_caracteriza} to obtain isolating blocks $N$ which are topological $2$--manifolds) and $K$ has the shape of $M$ with $r$ nonseparating circumferences removed (recall that $n^{-}$ does not separate $M$, as shown is Step 2 in the proof of Theorem \ref{teorema_cohomologia}). This has the same homotopy type as a compact connected $2$--manifold whose boundary consists of $2r$ circumferences. Pasting $2r$ disks to get rid of the boundary one obtains a compact connected $2$--manifold without boundary, which has a standard representation as a polygon with edges identified. Now remove again the $2r$ disks from the interior of this polygon and observe that the remaining structure has the homotopy type of a wedge sum of circumferences. Therefore $K$ has the shape of a bouquet of, say, $q \geq 1$ circumferences.

Next we shall determine the value of $q$. Since $\mathcal{A}(K) = M$, by Corollary \ref{corolario_caracteristica} and the fact that the dimension of $M$ is even we have $\chi(M) = \chi(K)$. However $\chi(K) = 1 - q$ because $K$ is a wedge sum of $q$ circumferences and it follows that $q = 1 - \chi(K) = 1 - \chi(M)$. Finally we have $1 \leq q = 1 - \chi(M)$ which implies $\chi(M) \leq 0$.
 \end{proof}

It is not difficult to construct specific examples, for any compact surface $M$ without boundary such that $\chi(M) \leq 0$, of connected isolated unstable attractors with no external explosions which are global. Thus the result above is sharp.

Let us remark finally that the fact that neither the sphere nor the projective plane can contain connected attractors without external explosions will also be proved, by other means, in Examples \ref{ejemplo_esfera} and \ref{ejemplo_proyectivo}. 

\section{Manifolds which contain attractors with no external explosions}

Now we start our global study, relating dynamical properties of attractors with no external explosions (namely, the number of homoclinic components in their basin of attraction) with the topology of the phase space. An instance of this fact was presented in \cite[Theorem 17]{moronsanjurjoyo1}, where it was shown that every connected unstable attractor in $\mathbb{R}^n$ must have external explosions.

Let us remark here that, although this approach is inscribed in the lines of the classical work by Morse, Smale and Conley, it is not subsumed in it. Indeed, the Morse equation of a Morse decomposition $\mathcal{M}$ provides a relation between the Conley indices of the Morse sets in $\mathcal{M}$ and the Poincar\'e polynomial of $M$, but unstable attractors cannot be detected by Morse decompositions.

More precisely, suppose that $\mathcal{M} = \{M_1, \ldots, M_n\}$ is a Morse decomposition of $M$ and $K$ is an attractor contained in some $M_k$. Then we claim that $\widehat{K} \subseteq M_k$ also, which means that $K$ and $\widehat{K}$ are indistinguishable for $\mathcal{M}$. Indeed if $p \in \widehat{K}$, then $\omega(p) \subseteq K \subseteq M_k$ and, since $\alpha(p)$ must be contained in some Morse set and intersects $K$ (see the characterization of $\widehat{K}$ in \cite{bhatiaszego}) then $\alpha(p) \subseteq M_k$ also. This means that $p \in M_k$ too and proves that $\widehat{K} \subseteq M_k$, as claimed. In particular our results cannot be proved within the classical Morse--Conley theory.

\subsection{Necessary conditions for a manifold to contain an unstable attractor with no external explosions}

\begin{theorem} \label{teorema_cohomologia} Let $M$ be an $R$--orientable $n$--manifold without boundary. Assume that $K \subseteq M$ is a connected isolated unstable attractor with only internal explosions. Then there exist $r$ independent cohomology classes \[\alpha_1, \ldots, \alpha_r \in H^1(M;R)\] such that \[\alpha_i \smile \alpha_j = 0\ \forall\ 1 \leq i,j \leq r,\] where $r \geq 1$ is the number of homoclinic components of $\mathcal{A}(K) - K$.
\end{theorem}
\begin{proof} Since $K$ and thus $\mathcal{A}(K)$ are connected, we can assume without loss of generality that $M$ is itself connected (otherwise we can argue with the component which contains $K$), hence a metrizable locally compact space. Let $N = N^{+} \cup N^{-}$ be a connected regular isolating neighbourhood for $K$. The hypothesis that $K$ be unstable implies that $n^{-}$ is nonempty and has components $n_1^{-}, \ldots, n_r^{-}$ (the notation is the same as in Lemma \ref{describe_a_de_k}). Coefficients in $R$ are to be understood throughout.
\medskip

{\it Step 1.\/} $n^{-}$ does not separate $M$.
\begin{proof} It will be enough to show that it does not separate $\mathcal{A}(K)$ because the latter is an open neighbourhood of $n^{-}$. By Lemma \ref{describe_a_de_k} we have $\mathcal{A}(K) = K \cup n^{-} \cdot \mathbb{R} \cup \bigcup_{j=r+1}^s n_j^{+} \cdot \mathbb{R}$, where the union is disjoint. Hence \[\mathcal{A}(K) - n^{-} = K \cup n^{-} \cdot (-\infty,0) \cup n^{-} \cdot (0,+\infty) \cup \bigcup_{j=r+1}^s n_j^{+} \cdot \mathbb{R}.\]

For any $x \in n^{-} \cdot (-\infty,0)$ clearly $x \cdot (-\infty,0] \subseteq n^{-} \cdot (-\infty,0) \subseteq N^{-}$ so $\emptyset \neq \alpha(x) \subseteq K$. Therefore $\overline{x \cdot (-\infty,0]} = x \cdot (-\infty,0] \cup \alpha(x) \subseteq n^{-} \cdot (-\infty,0) \cup K \subseteq \mathcal{A}(K) - n^{-}$ is a connected subset of $\mathcal{A}(K) - n^{-}$ that contains $x$ and has nonempty intersection with $K$, which is also connected. This implies that $K$ and $x$ are in the same component of $\mathcal{A}(K) - n^{-}$. Similar arguments can be made for $x$ in $n^{-} \cdot (0,+\infty)$ or $\bigcup_{j=r+1}^s n_j^{+} \cdot \mathbb{R}$, the only difference being that now it is the forward semitrajectory $x \cdot [0,+\infty)$ and the $\omega$--limit $\emptyset \neq \omega(x) \subseteq K$ which are to be used. This proves that every point in $\mathcal{A}(K) - n^{-}$ lies in the same component as $K$, so that $\mathcal{A}(K) - n^{-}$ is indeed connected.
\end{proof}
\medskip

Now we shall define the cohomology classes $\alpha_1, \ldots, \alpha_r \in H^1(M)$. Let us denote $k : M \hookrightarrow (M,M-n^{-})$ the inclusion. Since $H_1(M,M-n^{-}) = H^{n-1}(n^{-})$ by Alexander duality on $M$ and the latter module is, by Lemma \ref{lema_dualidad}(2), isomorphic to $H_0(n^{-}) = R^r$, it follows that $H_1(M,M-n^{-})$ is free of rank $r$. Now by the universal coefficients theorem $H^1(M,M-n^{-})$ is also free of rank $r$. Let $\beta_1, \ldots, \beta_r \in H^1(M,M-n^{-})$ be independent generators and set $\alpha_j = k^{*}(\beta_j) \in H^1(M)$ for every $1 \leq j \leq r$.
\medskip

{\it Step 2.\/} The classes $\alpha_1, \ldots, \alpha_r \in H^1(M)$ are independent.
\begin{proof} Consider the following portion of the long exact sequence in reduced cohomology for the pair $(M,M-n^{-})$: \[\ldots \longleftarrow H^1(M) \stackrel{k^{*}}{\longleftarrow} H^1(M,M-n^{-}) \longleftarrow \widetilde{H}^0(M - n^{-}) \longleftarrow \ldots\] Since we know by Step 1 that $n^{-}$ does not separate $M$, $\widetilde{H}^0(M - n^{-}) = 0$ and consequently $k^{*} : H^1(M,M-n^{-}) \longrightarrow H^1(M)$ is a monomorphism. Therefore $\alpha_1, \ldots, \alpha_r$ are independent classes in $H^1(M)$.
\end{proof}
\medskip

{\it Step 3.\/} The classes $\alpha_1, \ldots, \alpha_r$ satisfy the relations $\alpha_i \smile \alpha_j = 0$ for every $1 \leq i,j \leq r$.
\begin{proof} Observe that $n^{-}$ possesses a product neighbourhood in $M$, namely $\bigcup_{j=1}^r C_j = n^{-} \cdot \mathbb{R} \cong n^{-} \times \mathbb{R}$. Therefore there exist two homeomorphisms $h_1,h_2 : M \longrightarrow M$ such that: (1) both are homotopic to the identity ${\rm id}_M$, (2) $h_1(n^{-} \cdot (-1) ) = n^{-}$ and $h_2( n^{-} \cdot (+1) ) = n^{-}$. Moreover the cup product pairing of $H^1(M,M-n^{-} \cdot (-1) )$ and $H^1( M,M-n^{-} \cdot (+1) )$ is trivial because $M-n^{-} \cdot (-1)$ and $M - n^{-} \cdot (+1)$ are open sets whose union is the whole manifold $M$ (see \cite[p. 251]{spanier}).

Denoting $k_1 : M \hookrightarrow (M,M-n^{-} \cdot (-1))$ and $k_2 : M \hookrightarrow (M,M-n^{-} \cdot (+1))$ the inclusions, by the naturality of cup product (\cite[8., p. 251]{spanier}) there exists a commutative diagram

\begin{diagram}
H^1(M,M-n^{-}) & & H^1(M,M-n^{-}) & & \\
\dTo <{h_1^{*}} & & \dTo <{h_2^{*}} & & \\
H^1(M,M-n^{-} \cdot (-1) ) & \times & H^1 (M,M-n^{-} \cdot (+1) ) & \rTo^{\smile} & 0 = H^1(M,M) \\
\dTo <{k_1^{*}} & & \dTo <{k_2^{*}} & & \dTo \\
H^1(M) & \times & H^1(M) & \rTo^{\smile} & H^2(M)
\end{diagram}

It follows that $(h_1 \circ k_1)^{*} (\beta_i) \smile (h_2 \circ k_2)^{*} (\beta_j) = (k_1^{*} \circ h_1^{*})(\beta_i) \smile (k_2^{*} \circ h_2^{*})(\beta_j) = 0$, but since $h_1$ and $h_2$ were chosen to be homotopic to the identity, we have $h_1 \circ k_1 \simeq k$, $h_2 \circ k_2 \simeq k$ and consequently $\alpha_i \smile \alpha_j = 0$.
\end{proof}
\medskip

This completes the proof.  \end{proof}

A number of useful corollaries follow inmediately from Theorem \ref{teorema_cohomologia}.

\begin{corollary} \label{corolario_inmediato} If $M$ is an $R$--orientable manifold such that $H^1(M;R) = 0$ then every connected isolated unstable attractor $K \subseteq M$ has external explosions.
\end{corollary}

As an instance of the result above we present the following Example, part of which was already contained in \cite[Theorem 17.]{moronsanjurjoyo1}. We would like to point out that, although the technique used there is quite different from ours, it has been inspiring for the present work.

\begin{example} \label{ejemplo_esfera} Every connected isolated unstable attractor $K \subseteq \mathbb{R}^n$ has external explosions, since $H^1(\mathbb{R}^n;\mathbb{Z}) = 0$. A similar statement holds for $\mathbb{S}^n$ (the $n$--sphere) when $n > 1$ or for the complex projective spaces $\mathbb{CP}^n$.
\end{example}

We have just concluded that every connected unstable attractor $K \subseteq \mathbb{R}^n$ must have external explosions, because $\mathbb{R}^n$ does not fulfill the necessary conditions given by Theorem \ref{teorema_cohomologia}. The following Example presents a technique which allows to show that some phase spaces share the same property, although they {\it do} fulfill the requirements of Theorem \ref{teorema_cohomologia}.

\begin{example} \label{ejemplo_anillo} Every connected isolated unstable attractor $K$ contained in the open $2$--dimensional annulus $\mathbb{A} = \{ x \in \mathbb{R}^2 : 1 < \| x \| < 2 \}$ has external explosions.

To prove this observe that $\mathbb{A}$ is embedded in $\mathbb{R}^2$, where every connected isolated unstable attractor has external explosions. Our plan is to extend the flow in $\mathbb{A}$ to the whole of $\mathbb{R}^2$.

Assume that $K$ has no external explosions and let $\widehat{K} \subseteq \mathbb{A}$ be the stabilization of $K$. Let $P$ be a compact, positively invariant neighbourhood of $\widehat{K}$ in $\mathbb{A}$. By a theorem of Beck (\cite[3. Theorem, p. 99]{beck}) we can modify the flow in the annulus making every point in $\mathbb{A} - {\rm int}(P)$ fixed and leaving ${\rm int}(P)$ unchanged, save an appropriate reparametrization. Under this new flow $K$ is still a connected isolated unstable attractor whose basin of attraction is ${\rm int}(P)$ and which does not have external explosions. Extending the flow to the whole $\mathbb{R}^2$ by letting every point outside $\mathbb{A}$ be fixed we see that $K$ is a connected isolated unstable attractor in $\mathbb{R}^2$ which does not have external explosions, contradicting Example \ref{ejemplo_esfera}.
\end{example}

This stands in contrast with the situation for the closed annulus, which certainly contains an unstable attractor with no external explosions (apply the construction described in Example \ref{ejemplo_general} with $Z = [0,1]$, the unit interval).

The argument presented in the Example above can be used in general to obtain the following refinement of Corollary \ref{corolario_inmediato}.

\begin{corollary} \label{corolario_refinado} Let $M$ be an $n$--manifold without boundary which can be embedded in an $R$--orientable $n$--manifold $\widetilde{M}$ such that $H^1(\widetilde{M};R) = 0$. Then every connected isolated unstable attractor $K \subseteq M$ has external explosions.
\end{corollary}

So far we have not made use of the multiplicative structure in the cohomology ring of $M$. An example which takes it into account is gained with the real projective spaces.

\begin{example} \label{ejemplo_proyectivo} A connected isolated unstable attractor $K$ in $\mathbb{RP}^n$, the $n$--dimensional real projective space, has external explosions whenever $n \geq 2$.

To prove this let $R = \mathbb{Z}_2$ in Theorem \ref{teorema_cohomologia}. It is known that \[H^{*}(\mathbb{RP}^n;\mathbb{Z}_2) = \frac{\mathbb{Z}_2[\alpha]}{\alpha^{n+1}}\] is the truncated polynomial algebra in one indeterminate $\alpha$. Therefore the only nontrivial class in $H^1(\mathbb{RP}^n;\mathbb{Z}_2)$ is $\alpha$, but its square $\alpha^2$ is nontrivial in $H^2(\mathbb{RP}^n;\mathbb{Z}_2)$ (recall we assumed that $n \geq 2$) so the conditions in Theorem \ref{teorema_cohomologia} cannot be met.

Observe that in case $n=1$ the space $\mathbb{RP}^1$ is homeomorphic to $\mathbb{S}^1$, which certainly admits isolated unstable attractors with no external explosions (let $K$ be a single point and the rest of $\mathbb{S}^1$ a homoclinic orbit).
\end{example}

\begin{example} \label{ejemplo_toro} Let $\mathbb{T}^n = \mathbb{S}^1 \times \ldots \times \mathbb{S}^1$ be the $n$--dimensional torus. If $K \subseteq \mathbb{T}^n$ is an isolated unstable attractor with no external explosions, then there exists exactly one homoclinic component in $\mathcal{A}(K) - K$.

The assertion follows from the structure of the cohomology ring $H^{*}(\mathbb{T}^n;\mathbb{Z})$, namely that of an exterior algebra with $n$ generators $\omega_1, \ldots, \omega_n \in H^1(\mathbb{T}^n;\mathbb{Z})$. In fact, if $\alpha_1 = \sum_{i=1}^n k_i \omega_i$ is a nontrivial class in $H^1(\mathbb{T}^n;\mathbb{Z})$ and $\alpha_2 = \sum_{j=1}^n l_j \omega_j$ is another class satisfying $\alpha_1 \smile \alpha_2 = 0$, then an easy computation shows $\alpha_1 \smile \alpha_2 = \sum_{1 \leq i < j \leq n} (k_i l_j - k_j l_i) (\omega_i \smile \omega_j) = 0$ which implies $k_i l_j - k_j l_i = 0$ for every $1 \leq i < j \leq n$. But these are exactly the $2$--minors of the matrix \[\left( \begin{array}{cccc} k_1 & k_2 & \ldots & k_n \\ l_1 & l_2 & \ldots & l_n \end{array} \right)\] so it has rank one and $\alpha_1, \alpha_2$ are linearly dependent. Hence the maximum number of classes in $H^1(\mathbb{T}^n;\mathbb{Z})$ which meet the conditions of Theorem \ref{teorema_cohomologia} is one and the assertion is proved.
\end{example}

The next result generalizes \cite[Theorem 18]{moronsanjurjoyo1} and is a joint consequence of Theorem \ref{teorema_caracteriza} and Theorem \ref{teorema_cohomologia}. We shall state it for attractors in $\mathbb{R}^2$, the euclidean plane, but it holds more generally for any $2$--manifold where unstable attractors must have external explosions (such as $\mathbb{S}^2$, $\mathbb{RP}^2$ or open annuli, for example).

\begin{corollary} \label{corolario_cuandoestable} Let $K \subseteq \mathbb{R}^2$ be a connected isolated attractor. Then, \[K \text{ is stable } \Leftrightarrow \chi(K) = \chi(\mathcal{A}(K)).\]
\end{corollary}
\begin{proof} The only nontrivial implication is $\Leftarrow$. Thus, assume $\chi(K) = \chi(\mathcal{A}(K))$. Then by Theorem \ref{teorema_caracteriza} $K$ does not have external explosions, so using Example \ref{ejemplo_esfera} it follows that $K$ cannot be unstable.
\end{proof}

Thus for example a connected, isolated global attractor in $\mathbb{R}^2$ is stable if, and only if, it has trivial shape.

\subsection{Sufficient conditions for a manifold to contain an unstable attractor with no external explosions}

To finish this section we shall present a partial converse to Theorem \ref{teorema_cohomologia}. Basically, it is an elaborated setting of Example \ref{ejemplo_general}.

\begin{theorem} \label{proposicion_construye} Let $M$ be a closed oriented smooth manifold. If $H^1(M;\mathbb{Z}) \neq 0$ then $M$ contains a connected isolated unstable global attractor having no external explosions.
\end{theorem}
\begin{proof} We assume without loss of generality that $M$ is connected and take coefficients in $\mathbb{Z}$. Let $\alpha \in H^1(M)$ be a nonzero cohomology class and denote by $z \in H_{n-1}(M)$ its (nonzero) Poincar\'e dual. By \cite[Th\'{e}or\`{e}me II.27., p. 55]{thom} there exists an oriented closed smooth hypersurface $Z \subseteq M$ such that, if $i : Z \hookrightarrow M$ denotes the inclusion, $i_{*}([Z]) = z$ where $[Z] \in H_{n-1}(Z)$ is a fundamental class of $Z$. If $Z$ is not connected, let it have components $Z_1, \ldots, Z_p$ with fundamental classes $[Z_1], \ldots, [Z_p]$ such that $[Z] = [Z_1] + \ldots + [Z_p]$. Since $i_{*}([Z]) = i_{*}([Z_1]) + \ldots + i_{*}([Z_p]) = z$ is nontrivial, some $i_{*}([Z_k])$ is nontrivial too (although not necessarily equal to $z$). We keep $Z_k$ (which we shall again call $Z$) and discard the remaining components.
\medskip

{\it Step 1.\/} $Z$ does not separate $M$.
\begin{proof} Consider the following portion of the long exact sequence in homology for the pair $(M,Z)$: \begin{multline*} \ldots \longrightarrow H_n(Z) = 0 \longrightarrow H_n(M) \longrightarrow H_n(M,Z) \longrightarrow \\ \longrightarrow H_{n-1}(Z) \stackrel{i_{*}}{\longrightarrow} H_{n-1}(M) \longrightarrow \ldots \end{multline*}

The inclusion induced homomorphism $i_* : H_{n-1}(Z) \longrightarrow H_{n-1}(M)$ is injective, because it carries the generator $[Z]$ of $H_{n-1}(Z) = \mathbb{Z}$ onto a nontrivial element in $H_{n-1}(M)$ by construction and the latter is torsionfree since $M$ is orientable. It follows that $H_n(M,Z) = H_n(M) = \mathbb{Z}$ and by the universal coefficients theorem the free part of $H^n(M,Z)$ is isomorphic to $\mathbb{Z}$.

From Alexander duality applied to the pair $(M,Z)$ in $M$ there exists an isomorphism $H_0(M-Z) = \check{H}^n(M,Z)$. Moreover, since $(M,Z)$ is a polyhedral pair, \v{C}ech and singular cohomology agree on it so $H_0(M-Z) = H^n(M,Z)$. Therefore $H^n(M,Z)$ is free, which together with the conclusion of the previous paragraph shows $H^n(M,Z) = \mathbb{Z}$. Hence $M-Z$ is connected.
\end{proof}
\medskip

{\it Step 2.\/} $Z$ possesses a product neighbourhood in $M$. More precisely, there exist an open neighbourhood $U$ of $Z$ in $M$ and a homeomorphism $h : U \longrightarrow Z \times \mathbb{R}$ such that $h(z) = (z,0)$ for every $z \in Z$.
\begin{proof} Consider the normal bundle of $Z$ in $M$, denoted by $\bot{Z}$. By the tubular neighbourhood theorem we can assume that its total space $E(\bot{Z})$ is embedded as an open neighbourhood $U$ of $Z$ in $M$ and $Z$ corresponds to the zero section of $\bot{Z}$. Thus it will be enough to show that $\bot{Z}$ is trivial. Now $\bot{Z}$ is isomorphic to the quotient bundle $\frac{TM|_{Z}}{TZ}$, where $TM$ and $TZ$ are the tangent bundles to $M$ and $Z$ respectively and $TM|_{Z}$ is the restriction of $TM$ to $Z$ (more precisely, the pullback of $TM$ under the inclusion $i : Z \hookrightarrow M$). Since $Z$ and $M$ are orientable, their tangent bundles are orientable and therefore so is $\bot{Z}$. Consequently $\bot{Z}$ is a one--dimensional orientable bundle, whence it is trivial (\cite[Theorem 4.3., p. 106]{hirsch}).
\end{proof}
\medskip

Finally, define a flow in $Z \times \mathbb{R}$ such that $Z \times (-\infty,0]$ and $Z \times [1,+\infty)$ consist of fixed points and points in $Z \times (0,1)$ move from $Z \times \{0\}$ to $Z \times \{1\}$. Carry this flow to $M$ via the homeomorphism $h$ and extend it letting every point outside $U$ be fixed. The set $K = M - h^{-1}(Z \times (0,1))$ is an isolated compact connected global attractor (it is connected because the inclusion $K \hookrightarrow M - Z$ is a shape equivalence by Lemma \ref{describe_a_de_k} and we have already shown that $Z$ does not separate $M$). $K$ is not stable, because the orbit of any $x \not\in K$ is homoclinic, and further it has no external explosions because by construction $J^{+}(x) \subseteq K$ for $x \in M - K$.
\end{proof}

Combining Theorems \ref{teorema_cohomologia} and \ref{proposicion_construye} the following characterization is inmediate.

\begin{corollary} Let $M$ be a closed orientable smooth manifold. There exists a connected isolated unstable attractor with no external explosions in $M$ if, and only if, $H^1(M;\mathbb{Z}) \neq 0$.
\end{corollary}

\section{Unstable attractors with the shape of $\mathbb{S}^n$ and no external explosions}

This last section abounds in the relation between properties of an attractor and those of the phase space. Namely, we tackle the following question: let $M$ be a connected manifold without boundary and let $K \subseteq M$ be a connected isolated unstable attractor with only internal explosions. If $K$ has the shape of $\mathbb{S}^n$, what can be said about $M$? An argument similar to that in Example \ref{nopuedeser} shows that $M$ must be $(n+1)$--dimensional. Further, by Corollary \ref{corolario_global} it follows that $M$ is compact and there exists just one homoclinic component. However now the problem is considerably more complicated because there does not exist a complete classification of higher--dimensional manifolds.

To state the result we need some notation. The product $\mathbb{S}^n \times \mathbb{S}^1$ can be thought of as the quotient obtained from $\mathbb{S}^n \times [-1,1]$ by pasting its upper and lower lids in standard fashion, that is identifying $(x,1)$ with $(x,-1)$. However we shall need a ``twisted product'' $\mathbb{S}^n \times_t \mathbb{S}^1$, which is the result of pasting the upper and lower lids of $\mathbb{S}^n \times [-1,1]$ via the orientation reversing homeomorphism $r_n : \mathbb{S}^n \longrightarrow \mathbb{S}^n$ given by $r_n(x_1, x_2, \ldots, x_{n+1}) = (-x_1, x_2, \ldots, x_{n+1})$, that is identifying $(x,1)$ with $(r_n(x),-1)$. Both $\mathbb{S}^n \times \mathbb{S}^1$ and $\mathbb{S}^n \times_t \mathbb{S}^1$ are compact, connected $(n+1)$--manifolds without boundary, but the first is orientable whereas the second is not.

We shall also assume that Poincar\'{e}'s conjecture is true throughout.

\begin{theorem} \label{teorema_esfera} Let $M$ be a connected manifold without boundary and let $K \subseteq M$ be an isolated unstable attractor without external explosions. If $K$ has the shape of $\mathbb{S}^n$, where $1 \leq n \leq 2$, then $M$ is homeomorphic to $\mathbb{S}^n \times \mathbb{S}^1$ (if orientable) or to $\mathbb{S}^n \times_t \mathbb{S}^1$ (if nonorientable) and $K$ is a global attractor.
\end{theorem}

Further generalizations (for arbitrary $n$) will be described later. To prove the Theorem we need the following Lemma \ref{lema_esfera}. Its proof relies upon the {\sl annulus theorem}, which states that if $f_1, f_2 : \mathbb{S}^{n-1} \longrightarrow \mathbb{S}^n$ are two bicollared embeddings of $\mathbb{S}^{n-1}$ into $\mathbb{S}^n$, then the component of $\mathbb{S}^n - f_1(\mathbb{S}^{n-1}) \cup f_2(\mathbb{S}^{n-1})$ comprised between $f_1(\mathbb{S}^{n-1})$ and $f_2(\mathbb{S}^{n-1})$ is a topological annulus, that is, a product $\mathbb{S}^{n-1} \times [0,1]$ (an embedding $f : \mathbb{S}^{n-1} \longrightarrow \mathbb{S}^n$ is bicollared if it can be extended to an embedding $\widehat{f} : \mathbb{S}^{n-1} \times [-1,1] \longrightarrow \mathbb{S}^n$ which restricts to $f$ on $\mathbb{S}^{n-1} \times \{0\}$). The annulus theorem was proved in \cite{kirby1} for $n \neq 4$ and in \cite{quinn1} for $n = 4$ (see also \cite{edwards1}).

\begin{lemma} \label{lema_esfera} Let $N$ be a compact $(n+1)$--manifold which has the shape of $\mathbb{S}^n$, $n \geq 2$. Then $N$ is orientable and has a boundary $\partial N$ with two components, which are homology $n$--spheres. If both are topological $n$--spheres, $N$ is homeomorphic to $\mathbb{S}^n \times [0,1]$.
\end{lemma}
\begin{proof} Observe in the first place that, since $N$ and $\mathbb{S}^n$ are ANR's, the hypothesis that both have the same shape implies that they have, in fact, the same homotopy type. In particular $N$ is simply connected (recall $n \geq 2$) so it is orientable. By Lefschetz duality (\cite[20 Theorem p. 298]{spanier}) $H_k(N,\partial N) = H^{n+1-k}(N) = H^{n+1-k}(\mathbb{S}^n)$ so $H_k(N,\partial N) = \mathbb{Z}$ for $k = 1$ or $k = n + 1$ and $H_k(N,\partial N) = 0$ otherwise. When $0 < k < n$ we have $H_k(N) = 0$ and $H_{k+1}(N,\partial N) = 0$ so from the exact sequence for the pair $(N,\partial N)$ \[\ldots \longrightarrow H_{k+1}(N,\partial N) = 0 \longrightarrow H_k(\partial N) \longrightarrow H_k(N) = 0 \longrightarrow \ldots\] it follows that $H_k(\partial N) = 0$. In the lower dimensions \begin{multline*} \ldots \longrightarrow H_1(N) = 0 \longrightarrow H_1(N,\partial N) = \mathbb{Z} \longrightarrow \\ \longrightarrow \widetilde{H}_0(\partial N) \longrightarrow \widetilde{H}_0(N) = 0 \longrightarrow \ldots \end{multline*} and we obtain $\widetilde{H}_0(\partial N) = \mathbb{Z}$, so $\partial N$ has two connected components $S_1$ and $S_2$. Since $H_k(\partial N) = H_k(S_1) \oplus H_k(S_2)$, for $0 < k < n$ we have $H_k(S_1) = H_k(S_2) = 0$. Finally, $\partial N$ is a compact, orientable, boundaryless $n$--manifold (being the boundary of a compact, orientable $(n+1)$--manifold), so the same is true for its components $S_1$ and $S_2$. Thus $H_n(S_1) = H_n(S_2) = \mathbb{Z}$ and consequently $S_1$ and $S_2$ are homology spheres.

Assume now that both homology $n$--spheres $S_1$ and $S_2$ are in fact $n$--spheres. Denote $\widehat{N}$ the boundaryless compact $(n+1)$--manifold obtained by pasting two disjoint $(n+1)$--balls $D_1$ and $D_2$ onto $N$ so that $\partial D_1 = S_1$ and $\partial D_2 = S_2$.

We claim that $\widehat{N}$ is a homology $(n+1)$--sphere. From the exact sequence for the pair $(\widehat{N},D_1 \cup D_2)$ it follows inmediately that $H_1(\widehat{N}) \oplus \mathbb{Z} = H_1(\widehat{N},D_1 \cup D_2)$ and $H_k(\widehat{N}) = H_k(\widehat{N},D_1 \cup D_2)$ for $2 \leq k \leq n+1$. Letting $p_1$ and $p_2$ denote the centers of $D_1$ and $D_2$ we have $H_k(\widehat{N},D_1 \cup D_2) = H_k(\widehat{N} - \{p_1,p_2\}, D_1 \cup D_2 - \{p_1,p_2\}) = H_k(N,\partial N) = H^{n+1-k}(N)$ by excision and Lefschetz duality. Therefore $H_1(\widehat{N}) \oplus \mathbb{Z} = H^n(N)$ and $H_k(\widehat{N}) = H^{n+1-k}(N)$ for $2 \leq k \leq n+1$. Since $N$ has the same cohomology groups of $\mathbb{S}^n$, substituting in the formulae above it is readily seen that $\widehat{N}$ has the homology groups of the $(n+1)$--sphere.

$\widehat{N}$ is obtained adjoining the $(n+1)$--cells $D_1$ and $D_2$ onto $N$. Since $n+1 \geq 3$ and $N$ is simply connected, so is $\widehat{N}$ (this is a consequence of Van Kampen's theorem \cite[Theorem 1.20.]{hatcher1}). Hence by the Poincar\'e conjecture $\widehat{N}$ is the $(n+1)$--sphere. Now observe that $S_1$ and $S_2$ are two embeddings of $\mathbb{S}^n$ into $\widehat{N} \cong \mathbb{S}^{n+1}$ which are bicollared. Namely, the union of a collar of $S_1$ in $N$ (see for example \cite[Proposition 3.42.]{hatcher1}) and another collar for $S_1$ in $D_1$ gives a bicollar for $S_1$ in $\widehat{N} \cong \mathbb{S}^{n+1}$. Similarly for $S_2$. From the annulus theorem it follows that $N$ is homeomorphic to $\mathbb{S}^n \times [0,1]$.
 \end{proof}

The Lemma just proved is not true for $n = 1$, because then $N$ can be either $\mathbb{S}^1 \times [0,1]$ or a closed M\"obius band (this is trivial).

\begin{proof} (of Theorem \ref{teorema_esfera}) We shall deal with the case $n = 2$, because $n = 1$ is similar but simpler. We already know that $M$ is a compact $3$--manifold and $K$ is a global attractor. Let $N$ be an isolating block for $K$ as in Lemma \ref{describe_a_de_k}.
\medskip

{\it Step 1.} $N$ is a $3$--manifold.
\begin{proof} Since $K$ is a global attractor and $M$ is compact, $K$ is purely unstable. Therefore both $n^{-}$ and $n^{+}$ are sections of $M - K$, so since the latter is a $3$--manifold and therefore a generalized $3$--manifold, it follows that $n^{-}$ and $n^{+}$ are generalized $2$--manifolds (as in Theorem \ref{teorema_superficies}). The same argument as in Theorem \ref{teorema_superficies} shows that $n^{-}$ and $n^{+}$ are compact, boundaryless, topological $2$--manifolds. The fact that $n^{-} \times (-1,0] \cong n^{-} \cdot (-1,0] \subseteq N$ is an open neighbourhood of $n^{-}$ in $N$ and similary for $n^{+}$ proves that $\partial N$ is properly attached to $N$, which is therefore a manifold.
\end{proof}
\medskip

{\it Step 2.} There exists a homeomorphism $h_1 : \mathbb{S}^2 \times [0,1] \longrightarrow N$ such that $h_1(\mathbb{S}^2 \times \{0\}) = n^{-}$ and $h_1(\mathbb{S}^2 \times \{1\}) = n^{+}$.
\begin{proof} Since the inclusion $K \hookrightarrow N$ is a shape equivalence, $N$ has the shape of $K$, hence that of $\mathbb{S}^n$. We know (by Lemma \ref{lema_esfera}) that the boundary components of $N$ are homology spheres, and the only $2$--manifold which satisfies this condition is the $2$--sphere. Hence again by Lemma \ref{lema_esfera} it follows that $N$ is homeomorphic to $\mathbb{S}^2 \times [0,1]$ via some $h_1 : \mathbb{S}^2 \times [0,1] \longrightarrow N$. Since $h_1$ restricts to a homeomorphism between the boundaries of $N$ and $\mathbb{S}^2 \times [0,1]$, there is clearly no loss of generality in assuming the required conditions $h_1(\mathbb{S}^2 \times \{0\}) = n^{-}$ and $h_1(\mathbb{S}^2 \times \{1\}) = n^{+}$ because if this does not happen it is only necessary to precede $h_1$ by $\mathbb{S}^2 \times [0,1] \ni (x,s) \mapsto (x,1-s) \in \mathbb{S}^2 \times [0,1]$.
\end{proof}
\medskip

{\it Step 3.} $M$ is homeomorphic either to $\mathbb{S}^2 \times \mathbb{S}^1$ or to $\mathbb{S}^2 \times_t \mathbb{S}^1$.
\begin{proof} It is easy to see that we can assume $n^{-} = n^{+} \cdot (-1)$ because both are sections of $M - K$. Denoting $P = n^{+} \cdot [-1,0]$ the part of $M$ comprised between $n^{-}$ and $n^{+}$, it is homeomorphic to $n^{-} \times [-1,0]$, hence to $\mathbb{S}^2 \times [-1,0]$, via a homeomorphism $h_2 : \mathbb{S}^2 \times [-1,0] \longrightarrow P$. As before we can be assume that $h_2(\mathbb{S}^2 \times \{-1\}) = n^{+}$ and $h_2(\mathbb{S}^2 \times \{0\}) = n^{-}$ without loss of generality.

Consider the compositions \[k : \mathbb{S}^2 \times \{0\} \stackrel{h_2}{\longrightarrow} n^{-} \stackrel{h_1^{-1}}{\longrightarrow} \mathbb{S}^2 \times \{0\}\] and \[\ell : \mathbb{S}^2 \times \{-1\} \stackrel{h_2}{\longrightarrow} n^{+} \stackrel{h_1^{-1}}{\longrightarrow} \mathbb{S}^2 \times \{1\},\] which can be thought of as self homeomorphisms of the $2$--sphere. Therefore either they are isotopic, or $k$ and $\ell \circ r_2$ are isotopic. This is a consequence of the {\sl stable homeomorphism conjecture} for the sphere, which asserts that every orientation preserving homeomorphism of the sphere is stable (hence isotopic to the identity, see for example \cite[Section 4.]{browngluck1}). The conjecture is now known to be true because it follows from the annulus conjecture, see \cite[Theorem 9.4.]{browngluck1}.

We shall assume that $k$ and $\ell \circ r_2$ are isotopic, for example (the other possibility is similar but somewhat simpler), so that there exists a level preserving homeomorphism $G : \mathbb{S}^2 \times [0,1] \longrightarrow \mathbb{S}^2 \times [0,1]$ such that $G_0(x) = G(x,0) = k(x,0)$ and $G_1(x) = G(x,1) = \ell(r_2(x),-1)$. Let $h_1^* = h_1 \circ G : \mathbb{S}^2 \times [0,1] \longrightarrow N$, which is a new homeomorphism with the same boundary properties as $h_1$, because $G$ preserves levels.

Now $h_1^*$ and $h_2$ agree on $\mathbb{S}^2 \times \{0\}$, because $h_1^*|_{\mathbb{S}^2 \times \{0\}} = h_1 \circ G_0 = h_1 \circ k = h_2$. Therefore both can be pasted to obtain a surjective continuous mapping $H : \mathbb{S}^2 \times [-1,1] \longrightarrow M$ defined by $H |_{\mathbb{S}^2 \times [-1,0]} = h_2$ and $H |_{\mathbb{S}^2 \times [0,1]} = h_1^{*}$. Observe further that $M$ is thus exhibited as a quotient of $\mathbb{S}^2 \times [-1,1]$ by identifying points which have the same image under $H$. These are precisely those $(x,1)$ and $(y,-1)$, of the upper and lower lids respectively, such that $H(x,1) = H(y,-1)$. However $H(x,1) = h_1^*(x,1) = h_1 \circ G_1(x,1) = h_1 \circ \ell \circ (r_2(x),-1) = h_2(r_2(x),-1)$ and $H(y,-1) = h_2(y,-1)$ so $H(x,1) = H(y,-1)$ if, and only if, $y = r_2(x)$ because $h_2$ is a homeomorphism. Therefore $M$ is the result of identifying the upper and lower lids of $\mathbb{S}^2 \times [-1,1]$ by means of $r_2$, so it is homeomorphic to $\mathbb{S}^2 \times_t \mathbb{S}^1$.
\end{proof}

Naturally, if $k$ and $\ell$ had been isotopic, then $M$ would have been homeomorphic to $\mathbb{S}^2 \times \mathbb{S}^1$.
 \end{proof}

Observe that Theorem \ref{teorema_esfera}, for the case $n=1$, can also follow from the results presented in the section about surfaces: if $K$ has the shape of (a wedge sum) a one--sphere, then the ambient manifold must be either the torus (which is the untwisted product $\mathbb{S}^1 \times \mathbb{S}^1$) or the Klein bottle (which is the twisted product $\mathbb{S}^1 \times_t \mathbb{S}^1$).

Now we shall describe the necessary modifications to obtain a higher--dimensional version of Theorem \ref{teorema_esfera}, that is when $n \geq 3$. The first problem one must handle is the fact that in Step 1 one cannot in general assure any more that $N$ is a manifold, because although its boundary is still a homology manifold it need not be a topological manifold anymore (see \cite[p. 245]{wilder}). This can be fixed assuming that $M$ is a differentiable manifold and the flow is smooth, since then the isolating block $N$ can be taken to be a manifold by \cite{easton}.

A second difficulty arises in the need to prove that $n^{-}$ and $n^{+}$, the components of $\partial N$, are simply connected. It suffices to guarantee that $\pi_1(n^{-})$ is abelian, since then it is trivial because $n^{-}$ is a homology sphere. The same happens for $n^{+}$ because it is homeomorphic to $n^{-}$, both being sections of $M - K$.

Thus

\begin{theorem} \label{teorema_esfera1} Let $M$ be a connected smooth manifold without boundary and let $K \subseteq M$ be an isolated unstable attractor with only internal explosions for a smooth flow. If $K$ has the shape of $\mathbb{S}^n$ and $M - K$ is simply connected (or more generally $\pi_1(M-K)$ is abelian), then $M$ is homeomorphic to $\mathbb{S}^n \times \mathbb{S}^1$ (if orientable) or to $\mathbb{S}^n \times_t \mathbb{S}^1$ (if nonorientable) and $K$ is a global attractor.
\end{theorem}

\section{Appendix} \label{apendice}

This appendix collects the proofs of two easy lemmata about algebraic topology which have been deferred to avoid disturbing the exposition.

\begin{lemma} \label{observacion_cociente} If $(X_1,A_1)$ and $(X_2,A_2)$ are compact pairs such that $X_1 - A_1 \cong X_2 - A_2$, then $p(X_1,A_1;R) = p(X_2,A_2;R)$.
\end{lemma}
\begin{proof} Denoting $\pi : (X_1,A_1) \longrightarrow (X_1/A_1,A_1/A_1)$ the canonical projection, by \cite[5. Theorem, p. 318]{spanier} we have that $\pi^{*}$ is an isomorphism so $p(X_1,A_1) = p(X_1/A_1,A_1/A_1)$, and by the same reason $p(X_2,A_2) = p(X_2/A_2,A_2/A_2)$. Clearly $X_1/A_1$ is a one--point compactification of $X_1 - A_1$ (the point at infinity being $A_1/A_1$) and similarly $X_2/A_2$ is a one--point compactification of $X_2 - A_2$, the point at infinity being $A_2/A_2$. However, since $X_1 - A_1 \cong X_2 - A_2$, the uniqueness of such a compactification implies that there exists a homeomorphism $h : (X_1/A_1,A_1/A_1) \longrightarrow (X_2/A_2,A_2/A_2)$. Consequently $p(X_1,A_1) = p(X_1/A_1,A_1/A_1) = p(X_2/A_2,A_2/A_2) = p(X_2,A_2)$.
 \end{proof}

\begin{lemma} \label{lema_suspension} For any topological space $X$, $H_k ( X \times [-1,1],X \times \{-1,1\};R ) = H_{k-1}(X;R)$. The same statement is true for cohomology.
\end{lemma}
\begin{proof} The proof is an application of Mayer--Vietoris sequences. Namely, consider the open pairs $( X \times [-1,1], X \times \{-1\} )$ and $( X \times [-1,1], X \times \{1\} )$ of subsets of the pair $( X \times [-1,1], X \times \{-1,1\} )$. Since the inclusions $X \times \{-1\} \hookrightarrow X \times [-1,1]$ and $X \times \{1\} \hookrightarrow X \times [-1,1]$ are homotopy equivalences, both pairs have zero homology modules in all dimensions. Moreover, their intersection is $(X \times [-1,1],\emptyset)$, which has homology isomorphic to that of $(X,\emptyset)$. Therefore their Mayer--Vietoris sequence reads \[\ldots \longrightarrow 0 \longrightarrow H_k(X \times [-1,1], X \times \{-1,1\} ) \longrightarrow H_{k-1}(X) \longrightarrow 0 \longrightarrow \ldots\] whence the result follows. For cohomology the same argument applies.
 \end{proof}

\bibliographystyle{amsplain}
\bibliography{unstable}
\end{document}